\documentclass[12pt,a4paper]{article}

\usepackage{amssymb,amsmath}
\usepackage{latexsym}
\usepackage[dvips]{graphicx}

\newcommand{\co}{{\mathbb C}}

\newcommand{\re}{{\mathbb R}}
\newcommand{\n}{{\mathbb N}}

\newcommand{\z}{{\mathbb Z}}
\newcommand{\cA}{{\cal{A}}}

\newcommand{\cB}{{\cal{B}}}

\newcommand{\cK}{{\cal{K}}}

\newcommand{\cS}{{\cal{S}}}

\newcommand{\cM}{{\cal{M}}}
\newcommand{\cI}{{\cal{I}}}

\newcommand{\fB}{{\mathbf B}}
\newcommand{\fC}{{\mathbf C}}
\newcommand{\fE}{{\mathbf E}}
\newcommand{\fG}{{\mathbf G}}

\newcommand{\fP}{{\mathbf P}}
\newcommand{\fSd}{{\mathbf S}^{d-1}}
\newcommand{\fQ}{{\mathbf Q}}

\newcommand{\bk}{{\mathbf k}}
\newcommand{\bn}{{\mathbf n}}

\newtheorem{theorem}{Theorem}
\newtheorem{prop}{Proposition}
\newtheorem{lemma}{Lemma}
\newtheorem{cor}{Corollary}
\newtheorem{remark}{Remark}
\newtheorem{ex}{Example}
\newtheorem{defi}{Definition}
\newtheorem{conj}{Conjecture}

\date{}
\author{V.\,Yu.~Protasov, A.\,S.~Voynov,
\thanks{Dept. of Mechanics and Mathematics, Moscow State University,
Vorobyovy Gory, 119992, Moscow, {e-mail: \tt\small
v-protassov@yandex.ru, an.voynov@gmail.com }}}
\title{Matrix semigroups with constant spectral radius
\thanks{
The first author is supported by the RFBR grants No
13-01-00642 and 14-01-00332, and by the grant of Dynasty foundation; the  second author is supported  by the RFBR grant No 14-01-00332 А, the grant of Dynasty foundation, the Simons grant, and the Science Schools grant No NS-3682.2014.1.}}

\begin{document}
\maketitle

\begin{abstract}
Multiplicative matrix semigroups with constant spectral radius (c.s.r.) are studied and applied to
several problems of algebra, combinatorics, functional equations, and dynamical systems.  We show that
all such semigroups are characterized by means of irreducible ones. Each irreducible c.s.r. semigroup defines
walks on Euclidean sphere, all its nonsingular elements are similar (in the same basis) to orthogonal.
We classify all nonnegative c.s.r. semigroups and arbitrary low-dimensional semigroups. For higher dimensions,
we describe five classes and leave an open problem on completeness of that list. The problem of algorithmic recognition
of c.s.r. property is proved to be polynomially solvable for irreducible semigroups and undecidable for reducible ones.

\medskip

\smallskip

\begin{flushright}
\noindent  \textbf{AMS 2010} {\em subject
classification: } 47D03,  15A30, 15A60, 15B36
\end{flushright}

\end{abstract}
\bigskip

\date{}

\medskip

\begin{center}
\textbf{1. Introduction}
\end{center}

We consider multiplicative closed semigroups of real $d\times d$ matrices. For a nonempty compact family~$\cA$ of  $\, d\times d$ matrices, we denote by $\cS_{\cA}$ the semigroup generated by~$\cA$ by multiplications and taking closure.

\begin{defi}\label{d5}
A matrix semigroup~$\cS$ has constant spectral radius (in short, $\cS$ is c.s.r.) if
the spectral radius of all its elements is the same and nonzero. A compact matrix family
$\cA$ is called c.s.r. if it generates a c.s.r. semigroup~$\cS_{\cA}$.
\end{defi}
The spectral radius $\rho(A)$ of a matrix~$A$  is the maximal modulus of its eigenvalues.
Clearly, in any c.s.r. semigroup the spectral radius of all matrices equals~$1$.
\begin{defi}\label{d7}
A family of matrices~$\cA$ is irreducible if there is no proper linear subspace of~$\re^d$
invariant for all matrices from~$\cA$.
\end{defi}
The semigroup $\cS_{\cA}$ generated by a family~$\cA$ is irreducible precisely when so is~$\cA$. Semigroups with constant
spectral radius have several important characteristic properties listed below. Some of them concern only irreducible semigroups.

\medskip

\noindent \textbf{1.} If $\cS$ is a semigroup with multiplicative spectral radius,
i.e., ${\rho(AB) = \rho(A)\rho(B), \ A, B\in \cS}$, then the semigroup
$\{[\rho(A)]^{-1}\, A \quad  |  \quad  A \in \cS \, , \  \rho(A)\ne 0\}$ has constant spectral radius. Thus, the study of semigroups
with multiplicative spectral radius is essentially reduced to c.s.r. semigroups. This observation was put to good use in~\cite{OR}. Moreover, for irreducible semigroups, the submultilicativity of the spectral radius ($\rho(AB) \le \rho(A)\rho(B), \ A, B\in \cS$) is equivalent to its multiplicativity~\cite[Theorem~2.1]{LR}. This property of semigroups can be relaxed further to the so-called Rota condition~\cite[Theorem 4.3]{OR}.

\medskip

\noindent \textbf{2.} All finite matrix semigroups that do not contain zero matrices are c.s.r. For integer matrices,
  the converse it also true: every irreducible c.s.r. semigroup of integer matrices is finite. We analyse this aspect
  in more detail in Subsection~9.1.

\medskip

\noindent \textbf{3.} For every  irreducible c.s.r. semigroup~$\cS$, there exists a norm in~$\re^d$
 such that the induced operator norm of all matrices from~$\cS$ is one~\cite[Theorem~4.7]{OR}.
 Since all matrix norms are equivalent, we obtain the following criterion:
\smallskip

 {\em An irreducible semigroup~$\cS$
 is c.s.r. if and only if there are two positive constants $C_1, C_2$ such that $C_1 \le \|A\| \le C_2$
 for all $A \in \cS$}.
\smallskip

 Thus, for irreducible semigroups, the c.s.r. property is equivalent to boundedness
 from above and from below.

\medskip

\noindent \textbf{4.}
Let $\cA$ be a compact family of matrices. The boundedness property (item 3) of the semigroup~$\cS_{\cA}$
means equal asymptotic behavior of
products $\Pi = A_{n}\cdots A_{1}$ for all possible $\, A_j \in \cA, \, j =1, \ldots , n, \, n \in \n$. Thus, {\em if $\cA$ is irreducible, then the semigroup $\cS_{\cA}$
is c.s.r. if and only if $\ C_1 \le \|\Pi\| \le C_2$ for all products~$\Pi$ of matrices from~$\cA$.}

For finite families~$\cA$, this property is better expressed in terms of their $p$-radii:
\begin{defi}\label{d8}
For a finite  family~$\cA= \{A_1, \ldots , A_m\}$ and for a given $p \in \re \cup \{\pm \infty\}$ the $p$-radius $\rho_p(\cA)$
is defined as follows:
\begin{equation}\label{prad}
\rho_p \ = \ \lim_{k \to \infty} \Bigl(m^{-p}\ \sum_{d_1, \ldots , d_k}\|A_{d_k}\cdots A_{d_1}\|^{\, p} \Bigr)^{\, \frac{1}{p k}}
\end{equation}
with the usual modifications in the cases $p =0$ and $p = \pm \infty$.
\end{defi}

If there is a zero product of matrices from~$\cA$, then for $p \le 0$,  we set $\rho_p = 0$. The $p$-radius is a non-decreasing function in~$p$. The $p$-radii are also called {\em joint spectral characteristics} of~$\cA$. For one matrix, i.e., in case $m=1\, , \, \cA = \{A\}$,  they are all coincide with the spectral radius of~$A$. For a family of matrices, they are always different, apart from the case when
the family $r^{-1}\, \cA \, = \, \{r^{-1}A_1,
\ldots , r^{-1}A_m\}$ is c.s.r. for some number~$r$:

{\em The $p$-radii of a family~$\cA$ coincide for all~$p$ if and only if the family $r^{-1}\, \cA$
is c.s.r. for some $r > 0$. In this case $\rho_p(\cA) = r$ for all $\ p \in \re \cup \{\pm \infty\}$. }

This is proved in the next section (Proposition~\ref{p3}), along with more details of the joint spectral characteristics.
\bigskip

Thus, the c.s.r. property of a family~$\cA$ means equal asymptotic behaviour of all products
of matrices from~$\cA$. Such families and semigroups play a special role in many applications, where
nonhomogeneous matrix products are used (see Section~9). The simplest example of the c.s.r. property is any
semigroup of orthogonal matrices. Another example is
a semigroup of stochastic matrices (such semigroups are, however, reducible). Other nontrivial examples were analyzed  in~\cite{LR, MS, OR, Po} along with many properties of c.s.r. semigroup.

In this paper, we begin with characterizing arbitrary c.s.r. semigroups by means of irreducible
semigroups (Theorem~\ref{th5}). Then we analyze irreducible c.s.r. semigroups.
Theorem~\ref{th10} shows that all nonsingular elements of an irreducible
c.s.r. semigroup are orthogonal matrices in some (common) basis. In particular, all c.s.r. semigroups of nonsingular matrices
are, up to a linear transform, subgroups of~$O(d)$. So, this case is simple. However, the problem of classifying all
c.s.r. semigroups (possibly containing singular matrices) is surprisingly nontrivial. We prove that an irreducible semigroup~$\cS$ is c.s.r. if and only if it defines walks on an ellipsoid, i.e., there exists a point~$x\in \re^d, \, x \ne 0$, such that its orbit $\{Ax\ | \   A \in \cS\}$ lies on  the surface of an ellipsoid (Theorem~\ref{th20}). This
theorem asserts only the existence of such a point~$x$ and does not give any recipe how to find it, along with the ellipsoid,  for a given semigroup~$\cS$. Note also that
 the ellipsoid in this theorem cannot be replaced by another convex body. In the special case, when all
 matrices are nonnegative, a complete classification of c.s.r. semigroups is obtained in Section~6.
For the general case,
we describe five classes of c.s.r. semigroups and leave an open problem that this classification is complete.
The affirmative answer is proved for low dimensions, in Section~8. Since for higher dimensions, the problem is unsolved,
it is natural to address the question of algorithmical recognition of c.s.r. semigroups. Given a finite set~$\cA$ of rational matrices,
one needs to decide whether the semigroup $\cS_{\cA}$ is c.s.r. The answer is rather curious: for irreducible~$\cA$,
 there is a polynomial time algorithm (we present it is Section~7), while for reducible~$\cA$ this problem is algorithmicaly udecidable, even if~$\cA$ is a pair of nonnegative matrices.

 Finally, in Section~9, we consider applications. Our results are applied to
 five different problems. In \S 9.1 we analyze finite matrix semigroups. We prove that
 there is a polynomial time algorithm that for a finite set of integer matrices, under some mild assumption, decides whether it generates a finite semigroup.  In \S 9.2 we consider linear switching systems
 of ODE and describe systems with equal asymptotic growth of trajectories for all switching laws. Then we apply c.s.r. semigroups to the study of fractal curves, including refinable functions and wavelets. The main result
 describes fractal curves with constant local regularity. Finally, in \S 9.5, we solve an open problem formulated by B.Reznick
(1990) on the asymptotics of the Euler binary partition function.
\smallskip

\enlargethispage{2\baselineskip}

Throughout the paper, we denote by $(\cdot \, , \, \cdot)$ the standard inner product in~$\re^d$, by $\|x\|_2 = \sqrt{(x,x)}$
the Euclidean norm, by $O(d)$ the orthogonal group, by $\mathbf{S}^{d-1}$ the unit Euclidean sphere in~$\re^d$.
We assume a basis in~$\re^d$ is fixed and identify a linear operator with the corresponding matrix.

\medskip

\newpage

\begin{center}
\textbf{2. The c.s.r. property and the joint spectral characteristics}
\end{center}
\bigskip

The joint spectral characteristics  of finite matrix sets such as Lyapunov exponents, joint and lower spectral radii, etc., have a rich history and numerous applications (see~\cite{BW, DL, GP, H} and references therein).
In view of Definition~\ref{d8}, all joint spectral characteristics form a one-parametric family
(the $p$-radii). Each of them indicates the asymptotic growth of the $L_p$-mean of norms for products
of matrices of length~$k$ as $k \to \infty$. The limit~(\ref{prad}) always exists and does not depend on the matrix norm~\cite{P1}. For ${p= 0}$, formula~(\ref{prad}) is modified
as $\, \rho_0 \, = \, \lim_{k \to \infty} \bigl(\prod\limits_{d_1, \ldots , d_k}\|A_{d_k}\cdots A_{d_1}\| \bigr)^{\, \frac{1}{k 2^k }}$,
this is the {\em Lyapunov exponent}. For $p = \pm \infty$ we obtain the
{\em joint spectral radius} and the {\em lower spectral radius} respectively:
\begin{equation}\label{jsr}
\rho_{+\infty}\ = \ \lim_{k \to \infty} \max_{d_1, \ldots , d_k}\|A_{d_k}\cdots A_{d_1}\|^{\, 1/k} \ ;
\qquad \rho_{-\infty}\ = \ \lim_{k \to \infty} \min_{d_1, \ldots , d_k}\|A_{d_k}\cdots A_{d_1}\|^{\, 1/k}\, .
\end{equation}
 The $p$-radius is a non-decreasing (typically, strictly increasing) function
in~$p$.  However, there are exceptions, when
all $p$-radii coincide. For example, in case of one matrix~$\cA = \{A\}$, we have $\rho_p = \rho(A)$ for all~$p$.
The following result describes all families with equal $p$-radii:
\begin{prop}\label{p3}
 We have $\rho_p(\cA) = r > 0$ for all $p \in \re\cup \{\pm \infty\}$,  if and only if the semigroup generated by the family~$r^{-1}\cA = \{r^{-1}A_1,
\ldots , r^{-1}A_m\}$ is c.s.r.
\end{prop}
{\tt Proof.} It suffices to consider the case~$r=1$. The following assertions are well-known:
\begin{equation}\label{twoineq1}
\rho_{-\infty} \ \le \  \min_{d_1, \ldots , d_k}\rho^{1/k}(A_{d_k}\cdots A_{d_1}) \ \le \
      \max_{d_1, \ldots , d_k}\rho^{1/k} (A_{d_k}\cdots A_{d_1})  \  \le \ \rho_{+\infty}\,  , \qquad k \in \n\, ,
\end{equation}
and
\begin{equation}\label{twoineq2}
\min_{d_1, \ldots , d_k}\rho^{1/k}(A_{d_k}\cdots A_{d_1}) \ \to \rho_{-\infty} \ ; \qquad
      \max_{d_1, \ldots , d_k}\rho^{1/k} (A_{d_k}\cdots A_{d_1})  \  \to \ \rho_{+\infty}\qquad \mbox{as} \quad k \ \to \ \infty\,
\end{equation}
(see, for instance,~\cite{GP}). If $\rho_{-\infty} = \rho_{+\infty} = 1$, then $\rho(A_{d_k}\cdots A_{d_1}) = 1$
for every product of matrices from $\cA$, hence $\cS_{\cA}$ is c.s.r. Conversely,
if $\cS_{\cA}$ is c.s.r., then $\rho(A_{d_k}\cdots A_{d_1}) = 1$ for each product. Substituting in~(\ref{twoineq2}) and taking limit as $k \to \infty$,
we obtain $\rho_{-\infty} = \rho_{+\infty} = 1$.

{\hfill $\Box$}
\smallskip

The notions of the joint and lower spectral radii are directly extended to arbitrary compact
families~$\cA$, with the same properties including~(\ref{twoineq1}) and~(\ref{twoineq2}).
In the sequel we use the short standard notation for them: $\rho_{+\infty} (\cA) = \rho(\cA)$
(the joint spectral radius) and $\rho_{-\infty}(\cA) = \check \rho(\cA)$ (the lower spectral radius).
Note that they are well-defined only for compact families~$\cA$.
\medskip

\medskip

\begin{center}
\textbf{3. Arbitrary c.s.r. semigroups vs irreducible semigroups}
\end{center}
\bigskip

Most of known results on semigroups with constant and with multiplicative spectral radius deal with irreducible semigroups~\cite{LR, MS, OR, Po} .
In this section we show that the general case  can be characterized by the irreducible one.
The main result, Theorem~\ref{th5}, asserts that any reducible c.s.r. semigroup can be factored to
several irreducible semigroups of smaller dimensions, one of which is c.s.r. and the others are contractions.

Let us begin with several auxiliary results. We say that a matrix semigroup~$\cS$ has a {\em bounded spectrum} if there is a constant $C$
such that $\rho(A) \le C$ for all $A \in \cA$.

\begin{prop}\label{p5}
An irreducible matrix semigroup has a  bounded spectrum if and only it is bounded.
\end{prop}
The proof is in Appendix. This result generalizes \cite[Theorem 4.1]{OR},
where it was proved that any irreducible c.s.r. semigroup is bounded.
Note that the irreducibility assumption is essential, the corresponding example
is the semigroup generated by one matrix $A$ such that $(A)_{21} = 0$ and all other entries are ones.
An important consequence is that the joint spectral radius~$\rho(\cS)$ is well-defined for
any irreducible closed semigroup~$\cS$ with bounded spectrum. Moreover, $\rho(\cS) \le 1$, since all
products of matrices from~$\cS$ are bounded.
\begin{prop}\label{p7}
For any irreducible matrix semigroup~$\cS$ with a bounded spectrum, there is a norm~$\|\cdot \|$ in~$\re^d$
such that, in the induced operator norm, we have $\|A\|\le 1, \, A \in \cS$.
\end{prop}
{\tt Proof.} By Proposition~\ref{p5}, the semigroup $\cS$ is bounded, hence $\rho(\cS)$ is well-defined and does not exceed one.
Since~$\cS$ is irreducible, there exists Barabanov's norm for~$\cS$ in $\re^d$, for which
$\|A\| \le \rho(\cS)\, , \ A \in \cS$ (see~\cite{B1}). This completes the proof.

{\hfill $\Box$}
\medskip

For any reducible semigroup~$\cS$, there exists a basis in~$\re^d$ in which
   every matrix $A \in \cS$ has the following block upper triangular form:
 \begin{equation}\label{blocks}
A \quad = \quad \left(
\begin{array}{cccccc}
A^{(1)} & * &  \ldots &  * \\
0 & A^{(2)}&  * & \vdots \\
\vdots & {} &  \ddots  & * \\
0 &  \ldots &  0 & A^{(s)}
\end{array}
\right)\
\end{equation}
with square diagonal blocks of sizes~$d^{(i)}, \, i = 1, \ldots , s, \  \sum_{i=1}^s d^{(i)} = d$.
The locations of diagonal blocks and their sizes are the same for all matrices~$A \in \cS$.
For each $i=1, \ldots , s$, the semigroup $\cS^{(i)}$ formed by $\, d^{(i)}\times d^{(i)}$ matrices of the $i$th block is
irreducible. If $\cS$ is irreducible, we set~$s=1, d^{(1)} = d, \cS^{(1)} = \cS$.

For every family of matrices~$\cA$, not necessarily a semigroup, we use the same factorization~(\ref{blocks})
with  irreducible families~$\cA^{(i)}$ of matrices in the $i$th  diagonal
blocks. Clearly, $\rho(A) = \max_{i = 1, \ldots , s} \rho(A_i)$. If all families~$\cA^{(i)}$ are compact
(and hence, the joint spectral radius is well-defined), then $\rho(\cA) = \max_{i = 1, \ldots , s} \rho(\cA_i)$
(see~\cite{DL}).
 \begin{theorem}\label{th5}
 A matrix semigroup $\cS$ is c.s.r. if and only if, in factorization~(\ref{blocks}),
 one of the irreducible semigroups~$\cS^{(j)}$ is c.s.r., and
 the others satisfy~$\, \rho(\cS^{(i)}) \le 1$.
 \end{theorem}
 \begin{remark}\label{r5}
 {\em Since all semigroups~$\cS_i$ are irreducible and have bounded spectra
 (because $\rho(A^{(i)}) \le \rho(A) = 1\, , \ A \in \cS$), Proposition~\ref{p5} implies
 that they are all compact. Therefore, $\rho(\cS^{(i)})$ is well-defined.

 For a semigroup $\cS_{\cA}$ generated by a matrix family~$\cA$,  Theorem~\ref{th5} is formulated in the same way:
 $\cA$ generates a c.s.r. semigroup if and only if, in factorization~(\ref{blocks}),
 one of the irreducible families~$\cA^{(j)}$ generates a c.s.r. semigroup, and
 the others satisfy~$\, \rho(\cA^{(i)}) \le 1$.
  }
 \end{remark}
 \smallskip

{\tt Proof of Theorem~\ref{th5}.} Sufficiency is clear: if $\cS$ is factored to the form~(\ref{blocks})
with a c.s.r. semigroup in the $j$th diagonal block,
then, for every $A \in \cS$ we have $\rho (A) \, = \, \max_{i=1, \ldots , s} \rho (A_i)\, = \, \rho(A_j) \, = \, 1$,
and hence $\cS$ is c.s.r.
Indeed, $\rho(A_j)  =  1$, since $\cS^{(j)}$ is c.s.r., and $\, \rho (A_i) \le \rho(\cS^{(i)}) \le 1$, for $i \ne j$.

Necessity. If $\cS$ is c.s.r., then all semigroups~$\cS^{(i)}$
have bounded spectra ($\rho(A^{(i)}) \le \rho(A) = 1\, , \ A \in \cS$), and hence, by Proposition~\ref{p7},
the norms of all matrices of these semigroups are bounded uniformly by some constant~$C$.
 If some semigroup $\cS_j$ is c.s.r., then the proof is completed.
If, otherwise,
for every $i = 1, \ldots , s$, there is a matrix $A_i \in \cS$, whose $i$th block
$A^{(i)}_i$ has its spectral radius smaller than~$1$, then, for an arbitrary~$n$, we consider
the product $\Pi = A_1^{n}\cdots A_s^{n} \, \in \cS$. If $n$ is large enough, then,
$\|(A_i^{(i)})^n\|\, < \, 1/C^2$, for every $i= 1, \ldots , s$. Therefore,
for the $i$th block of $\Pi$, we have
$$
 \bigl\|\Pi^{(i)} \bigr\| \ = \  \bigl\|(A^{(i)}_1)^{n}\cdots (A^{(i)}_s)^{n}\bigr\| \  = \
\bigl\|(A^{(i)}_1)^{n}\cdots (A^{(i)}_{i-1})^{n}\bigr\| \cdot
   \bigl\|(A^{(i)}_i)^{n}\bigr\| \cdot   \bigl\|(A^{(i)}_{i+1})^{n}\cdots (A^{(i)}_s)^{n}\bigr\|
$$
$$
     \le \  C \, \bigl\|(A^{(i)}_i)^{n}\bigr\| \, C \ < \  1 \, .
$$
Thus, $\rho(\Pi^{(i)}) \, \le \, \|\Pi^{(i)}\|\, < \, 1$,  for all~$i$, and so $\, \rho(\Pi) \, < \, 1$, which contradicts
to the c.s.r. property of~$\cS$.

{\hfill $\Box$}
\medskip

\begin{ex}\label{ex10}
{\em If a semigroup~$\cS$ consists of row-stochastic matrices, then it obviously c.s.r. In this case
one can take $e = (1, \cdots , 1)^T \in \re^d$ as the first basis vector, and the other $d-1$ basis vectors span the
subspace~$e^{\, \perp}$. In this basis, all matrices from~$\cS$ get the form~(\ref{blocks})
with $d^{(1)}=1, d^{(2)} = d-1$, and for every $A \in \cS$ the one-dimensional matrix~$A^{(1)}$ is~$1$.
}
\end{ex}
\bigskip

\begin{center}
\textbf{4. Nonsingular elements in c.s.r. semigroups}
\end{center}

\bigskip

The set of nonsingular matrices in a c.s.r. semigroups has a simple structure:
\begin{theorem}\label{th10}
If an irreducible  matrix semigroup~$\cS$ has constant spectral radius, then
in a suitable basis in~$\re^d$, all nonsingular elements of~$\cS$ are orthogonal matrices.
\end{theorem}
Thus, if $\cS$ is c.s.r., then there is a basis in~$\re^d$ in which all nonsingular matrices of~$\cS$ are orthogonal.
Note that it is easily shown that each nonsingular matrix in a c.s.r. semigroup is similar to an orthogonal one.
What is nontrivial in Theorem~\ref{th10} that all those similarities are realized by the same linear transform.
\smallskip

\noindent {\tt Proof.} By Proposition~\ref{p7}, there is norm $\|\cdot\|$ in~$\re^d$ such that
$\|A\| \le 1, \ A \in \cS$, in the induced operator norm.  For the unit ball~$\fB\subset \re^d$ of that norm, we have  $A\, \fB \subset \fB\, , A\in \cS$.
Furthermore, it was shown in~\cite[Theorem 2.5]{OR} that each element of an irreducible c.s.r. semigroup is a direct sum of a nilpotent operator and an operator similar to orthogonal. Hence, every nonsingular
element of~$\cS$ is similar to orthogonal. Therefore, its determinant is one, and hence this operator preserves the volume.
Consequently, for every nonsingular matrix~$A \in \cS$, we have $A\, \fB \, = \, \fB$. Let now~$\fE$ be the John ellipsoid of the
convex body~$\fB$, i.e., an ellipsoid of the maximal volume contained in~$\fB$. From the uniqueness of the John ellipsoid~\cite{Z}
it follows that~$A\, \fE = \fE$. Taking the basis, in which $\fE$ is a Euclidean ball, we obtain that every nonsingular element
 $A \in \cS$ is orthogonal.

{\hfill $\Box$}
\smallskip

\begin{cor}\label{cor10}
Suppose a semigroup~$\cS$ is generated
by an irreducible  set of nonsingular matrices; then
$\cS$ is c.s.r. if and only if there is a basis in~$\re^d$ in which all elements of~$\cS$ are orthogonal
matrices.
 \end{cor}
 \begin{remark}\label{r20}
 {\em A semigroup is generated by nonsingular matrices may contain singular elements. They
 may appear after taking closure}.
 \end{remark}
{\tt Proof of Corollary~\ref{cor10}.} By Theorem~\ref{th10} there is a basis in which  all elements of the generating set are orthogonal.
Hence, all their products are also orthogonal, and all limit points of these products also are.

{\hfill $\Box$}

Since a c.s.r. semigroup of nonsingular matrices is actually a group~\cite{OR}, we obtain
\begin{cor}\label{cor15}
If a c.s.r. semigroup consists of nonsingular matrices, it is similar (in a common basis in~$\re^d$)
to a subgroup of~$O(d)$.
 \end{cor}
\smallskip

\begin{center}
\textbf{5. Irreducible c.s.r. semigroups define walks on Euclidean sphere}
\end{center}

\bigskip

The results of previous section may make a wrong impression that all c.s.r. semigroups
are easily classified. In fact, there is a rich variety of classes of such semigroups, the problem of their complete
classification is still unsolved (see Section~8). In view of Theorem~\ref{th10}, this is caused by singular matrices
in c.s.r. semigroups. The following criterion is actually a weaker version of Theorem~\ref{th10}
that holds in general case, without the non-singularity assumption.
\begin{theorem}\label{th20}
An irreducible semigroup~$\cS$ is c.s.r. precisely when there is a suitable basis and a
point $x \in \re^d$ such that $\|x\|_2 = 1$ and $\|Ax\|_2 = 1$, for all $A \in \cS$.
\end{theorem}
{\tt Proof.}  {\em Sufficiency.} By irreducibility, the set $\{Ax \ | \ A \in \cS\}$ contains
linearly independent vectors~$\{x_i\}_{i=1}^d$. Since $\|Ax_i\|_2 = 1$ for every $A \in \cS$, it follows that
in the basis $\{x_i\}_{i=1}^d$ all the matrices from~$\cS$ are uniformly bounded, hence $\rho(\cS) \le 1$.
On the other hand, if $\rho(A) < 1$ for some $A \in \cS$, then $\|A^nx\|_2 < 1$, for sufficiently large~$n$,
which contradicts the assumption. Thus, $\rho(A) = 1$ for all $A \in \cS$.

{\em Necessity.} We prove under the assumption that~$\cS$ is finitely generated; the general case then follow by the
standard compactness argument. Thus, assume $\cS$ is generated by a set $\cA = \{A_1, \ldots , A_m\}$.
Since $\cS$ is c.s.r. it follows that $\rho_1(\cA) = \rho_2(\cA) = 1$. For an irreducible set
$\cA$, there exists an ellipsoidal norm $\|\cdot\|$ such that $\bigl( \frac{1}{m}\, \sum_{i=1}^m \|A_iz\|^2\bigr)^{1/2}
\, = \, \rho_2(\cA)\, \|z\|$ for every $z \in \re^d$~\cite{P1}. In a suitable basis, this norm becomes Euclidean.
Since in our case $\rho_2 = 1$, we have $\bigl( \frac{1}{m}\, \sum_{i=1}^m \|A_iz\|_2^2\bigr)^{1/2}
\, = \, \|z\|_2\, , \  z \in \re^d$. Iterating $k$ times, we get
 $\bigl( \frac{1}{m^k}\, \sum_{i_1, \ldots , i_k} \|A_{i_k}\ldots A_{i_1}z\|_2^2\bigr)^{1/2}
\, = \, \|z\|_2$. Applying now the inequality between the quadratic and arithmetic mean we
obtain
$$
f_k(z) \ = \ \frac{1}{m^k}\, \sum_{i_1, \ldots , i_k} \bigl\|\, A_{i_k}\ldots A_{i_1}z\, \bigr\|_2
\ \le  \ \bigl\|z\bigr\|_2\ , \qquad z \in \re^d\, .
$$
 This inequality becomes equality only if all the norms $\|A_{i_k}\ldots A_{i_1}z\|_2$
are equal to~$\|z\|_2$. Otherwise this inequality is strict. If it is strict for all $z \in \fSd$,
then there is $\gamma < 1$ such that $f_k(z)
\, \le  \, \gamma \, \|z\|_2\, , \ z \in \re^d$,  which yields  $\rho_1(\cA) \le \gamma^{1/k} < 1$.
 This contradiction shows that for each natural $k$, there is $z_k \in \fSd$
 such that $\|A_{i_k}\ldots A_{i_1}z_k\|_2 = \|z_k\|_2 = 1$
for all $i_1, \ldots ,i_k \in \{1, \ldots , m\}$.
  Note that by submultiplicativity of operator norm the function $f_k(z)$ in non-increasing in~$k$.
  Therefore, $f_k(z_k) = 1$ implies that $f_n(z_k) = 1$ for all $n \le k$, and hence
 $\|A_{i_n}\ldots A_{i_1}z_k\|_2  = 1$
for all $i_1, \ldots ,i_n \in \{1, \ldots , m\}$. Consequently, for every limit point
$x$ of the sequence $\{z_k\}_{k \in \n}$  we have $\|Ax\|_2 = 1, \, A \in \cS$.

{\hfill $\Box$}

\smallskip

\begin{remark}\label{r45}
{\em Theorem~\ref{th20} asserts only the existence of a point $x$ such that
all its images $Ax \, , \ A \in \cA$, lie on the surface of some ellipsoid. Thus, the point~$x$ makes walks on that ellipsoid under the action of the semigroup~$\cS$. If all matrices of~$\cS$ are nonsingular, then all
 points $x \in \re^d \setminus \{0\}$ possess this property and, moreover, the ellipsoids are homothetic for all 
 points.   (Theorem~\ref{th10}). In general, however, this is not the case. In Remark~\ref{r40} in Section~7 we shall see that
 the quadratic forms of all those ellipsoids are points of a common invariant affine plane of the tensor squares of matrices from~$\cS$. This gives a way to find those ellipsoids explicitly, for a given finitely generated semigroup~$\cS_{\cA}$. As for the points~$x$, we do not have a method to find at least one of them. }
\end{remark}

\begin{remark}\label{r50}
{\em Theorem~\ref{th20} holds only for Euclidean sphere, not for any other sphere of Banach norm in~$\re^d$.
To see this, it suffices to consider the irreducible c.s.r. semigroup~$O(d)$.}
\end{remark}

\begin{cor}\label{cor30}
An irreducible semigroup $\cS$ is c.s.r. if and only if there is a compact set $\Omega \in \fSd$ such that
$\, \Omega \, = \, \cup_{\, A \in \cS} A\, \Omega$.
\end{cor}
{\tt Proof.} The sufficiency follows directly from Theorem~\ref{th20}. To prove the necessity
we consider the set $\Omega_0 = \{Ax \ |  A \in \cS\}$, where $x \in \fSd$ is a point such that
$Ax \in \fSd$ for all $A \in \cS$. Let  $\Omega_k$ be the closure of the set $\cup_{\, A \in \cS} A\, \Omega_{k-1}\, , \ k \in \n$. Then $\Omega_k \subset \Omega_{k-1}\, , \, k \in \n$. The set $\Omega = \cup_{k \ge 0}\Omega_k$ possesses the desired property.

{\hfill $\Box$}
\smallskip


\begin{center}
\textbf{6. Nonnegative c.s.r. semigroups}
\end{center}
\bigskip

If all matrices in a semigroup are nonnegative (have nonnegative entries), the c.s.r. property can be efficiently characterized. The irreducibility assumption in this case is relaxed to the {\em positive irreducibility} (Definition~\ref{d10} below).

We call a {\em coordinate subspace} a subspace of~$\re^d$ spanned by $p$ basis vectors, where $1 \le p\le  d-1$.
A nonnegative matrix~$A$ is called {\em positively reducible}, if it has an invariant coordinate subspace.
In this case there is  a permutation of the canonical basis, after which~$A$ gets a block upper-triangular form with two  blocks $p\times p$ and $(d-p)\times (d-p)$ on the diagonal, where~$p$ is the dimension of the common invariant coordinate subspace.
Otherwise, if $A$ has no common invariant coordinate subspaces, it is referred to as {\em positively irreducible}.
In this case, for every~$i$ and~$j$, there is a power $k \ge 1$ such that $(A^k)_{ij} > 0$.
\begin{defi}\label{d10}
A family~$\cA$ of nonnegative matrices is called positively reducible if all its matrices share a common invariant coordinate subspace.  Otherwise~$\cA$ is called  positively irreducible.
\end{defi}
 If $\cA$ is  positively reducible, then there is  a permutation of the canonical basis after which all matrices from~$\cA$ get the block upper-triangular form~(\ref{blocks}) with square diagonal blocks of sizes~$d^{(i)}, \, i = 1, \ldots , s, \  \sum_{i=1}^s d^{(i)} = d$. The locations and sizes of the diagonal blocks are the same for all matrices~$A \in \cA$,
the families  $\cA^{(i)}$ in the blocks are all positively
irreducible. For positively irreducible families we still use factorization~(\ref{blocks})
 with~$s=1, d^{(1)} = d, \cA^{(1)} = \cA$. If $\cA = \cS$ is a semigroup, then we denote by $\cS^{(i)}$
the corresponding semigroups in the blocks.
The following proposition and theorem almost repeat Proposition~\ref{p5} and  Theorem~\ref{th5}
respectively, but for nonnegative semigroups and for
positive irreducibility. Their proofs are literally the same, replacing Barabanov's
theorem on invariant norms for irreducible families by  analogous~\cite[Theorem 3]{GP} for positively irreducible
families.
\begin{prop}\label{p130}
A positively irreducible matrix semigroup has a  bounded spectrum if and only it is bounded.
\end{prop}

 \begin{theorem}\label{th30}
 A nonnegative matrix semigroup $\cS$ is c.s.r. if and only if, in factorization~(\ref{blocks}),
 one of the positively irreducible semigroups~$\cS^{(j)}$ is c.s.r., and
 the others satisfy~$\, \rho(\cS^{(i)}) \le 1$.
 \end{theorem}
The following criterion ensures the c.s.r. property of a positively irreducible semigroup.
 \begin{theorem}\label{th40}
 A positively irreducible  matrix semigroup $\cS$ is c.s.r. if and only if
 there is a proper affine subspace $V \subset \re^d$ such that $0 \notin V\, , \, V \cap \re^d_+ \ne \emptyset$
 and $A V \subset V\, , \ A \in \cS$.
 \end{theorem}
The proof is in Appendix. Thus, an irreducible semigroup is c.s.r. precisely when it has
an invariant affine plane (maybe one-point) intersecting the positive orthant and not passing through the origin.
The following proposition ensures that every such a plane intersects an interior of positive orthant by a
bounded set.
\begin{prop}\label{p40}
Under the assumptions of Theorem~\ref{th40}, the set $\, \fP \, = \, V \cap \re^d_+$ is bounded (i.e., is a polyhedron)
and intersects the interior of~$\re^d_+$.
\end{prop}
The proof is in Appendix. Theorem~\ref{th40} and Proposition~\ref{p40}  give a simple classification
of nonnegative c.s.r. semigroups. One takes an arbitrary affine plane $V$ of dimension from $0$ to $d-1$
that does not pass through the origin and intersects the positive orthant by a bounded set. Then every semigroup of nonnegative matrices respecting this plane is has constant spectral radius, and vice versa: every positively irreducible c.s.r. semigroup is obtained this way.
Note that the problem of classification of arbitrary c.s.r. semigroups (without the nonnegativity assumption) is
much more difficult (Section~8).

\begin{ex}\label{ex20}
{\em If $\cS$ is a semigroup of row-stochastic matrices, then $L = \fP = \{e\}$, where $e = (1, \ldots , 1)^T\in \re^d$.
So, in this case ${\rm dim}\, L = 0$, although $L$ is nonempty, and the polyhedron $\fP$
is one point. If $\cS$ is a semigroup of column-stochastic matrices, then $V \, = \, \{x \in \re^d \ | \ (e, x) = 1\}$ and
$\fP$ is a simplex with vertices $e_1, \ldots , e_d$. }
\end{ex}
 \begin{cor}\label{c35}
If a semigroup of nonnegative matrices is c.s.r., then it is reducible.
 \end{cor}
{\tt Proof}. If $\cS$ is irreducible, then it is  positively irreducible as well, and
 Theorem~\ref{th30} gives a common affine invariant subspace~$L$ of matrices from~$\cS$
such that $0 \notin L $.
If ${\rm dim}\, L \, \ge 1$, then the linear part of $L$ is a common invariant subspace for~$\cS$,
if ${\rm dim}\, L \, = \, 0$, then so is the linear span of $\{0\}\cup L$. The contradiction completes the proof.

 {\hfill $\Box$}
\smallskip

Let us recall that 
 Corollary~\ref{c35} deals with usual reducibility, not positive one. Thus, if nonnegative matrices generate a c.s.r. semigroup, then they share a common invariant linear subspace. Irreducible semigroups of nonnegative matrices are never c.s.r. For arbitrary semigroups, not necessarily nonnegative,
 this is not true. Examples can be found in~\cite{OR} and in our Section~8.

Observe that if several matrices have a common affine subspace~$V$, then so does every their convex combination.
Invoking Theorem~\ref{th40} we obtain: if a positively irreducible semigroup~$\cS$ is c.s.r., then
${\rm co }(\cS)$ also is. Theorem~\ref{th30} extends this result to arbitrary nonnegative matrix semigroup.
This proves the following
\begin{theorem}\label{th50}
If $\cS$ is a nonnegative c.s.r. semigroup, then its convex hull~${\rm co }(\cS)$ also is.
 \end{theorem}
 The polyhedron $\fP = V\cap \re^d_+$ is invariant with respect to each matrix $A \in {\rm co }(\cS)$,
 i.e., $A\, \fP \, \subset \, \fP$. Applying the Brower fixed point theorem, we see that there exists $z \in \fP$ such that
 $Az = z$. Since $\rho(A) = 1$, it follows that $z$ is a Perron-Frobenius eigenvector for~$A$.
 \begin{cor}\label{c40}
 Let $\cS$ be a positively irreducible c.s.r. and~$V$ be its common invariant affine subspace from
 Theorem~\ref{th30}. Then the polyhedron $\fP = V\cap \re^d_+$ contains a Perron-Frobenius eigenvector of
each  matrix~$A \in {\rm co }(\cS)$.
 \end{cor}
The following proposition gives a method to decide the c.r.s. property for a given nonnegative semigroup.
\begin{prop}\label{p50}
Suppose $\cS$ is a nonnegative matrix semigroup, $\bar A \in {\rm co}(\cS)$ is a positively irreducible matrix,
$\rho(\bar A) =1$, and $v$ is its Perron-Frobenius eigenvector; then~$\cS$ is c.s.r. if and only if its matrices have
a common invariant linear subspace  that contains all vectors~${v - Av\, , \ A \in \cS}$, and does not contain~$v$.
\end{prop}
{\tt Proof.} If there is a positively irreducible matrix~$\bar A \in {\rm co}(\cS)$,~then $\cS$ is positively
irreducible, the Perron-Frobenius eigenvector $v$ of $\bar A$ is positive and unique up to normalization.
If there is a common invariant linear subspace~$\tilde V$  that contains all vectors~$v - Av\, , \ A \in \cS$
and does not contain~$v$,  then $V = v + \tilde V$ is a common invariant affine subspace for~$\cS$, and
$0 \notin L$. Moreover, $v \, \in \, V \cap \re^d_+$, hence $V \cap \re^d_+ \ne \emptyset$. Theorem~\ref{th40} now yields
that~$\cS$ is c.s.r. Conversely, if~$\cS$ is c.s.r., then so is ${\rm co}(\cS)$ (Theorem~\ref{th50}),
hence~${\rm co}(\cS)$ possesses a common invariant affine subspace~$V$, and (Corollary~\ref{c40})
the set $\fP = V\cap \re^d_+$ contains the Perron-Frobenius eigenvector~$v$ of $\bar A$. Therefore,~$\tilde V$
is a common invariant linear subspace of~$\cS$ containing all vectors~$v - Av\, , \, A \in \cS$.

{\hfill $\Box$}

\smallskip

\begin{remark}\label{r30}
{\em  All the results of this section can be generalized from nonnegative matrices to matrices that share a common invariant cone~$K$. This is an arbitrary convex solid pointed cone $K \subset \re^d$. In this case, the positive reducibility of~$\cS$ means that all matrices from~$\cS$ share a common invariant subspace which is a linear  span of a face of~$K$. Theorems~\ref{th40} and~\ref{th50}, Propositions~\ref{p40} and~\ref{p50}, and Corollaries~\ref{c35}, \ref{c40} remain true, if we replace~$\re^d$ by~$K$, and nonnegative matrices by matrices leaving~$K$ invariant. In particular, Corollary~\ref{c35} implies
}
\end{remark}

\begin{prop}\label{p60}
None of irreducible c.s.r. semigroups possess a common invariant cone.
\end{prop}

Let us remark  that for irreducible matrix semigroups, the equality $\rho_1(\cS) = \rho(\bar A)$,
where $\bar A = \frac1m \sum_{i=1}^m A_i\, ,$ and $A_1, \ldots , A_m$ are matrices generating~$\cS$, is necessary and sufficient for the existence  of a common invariant cone for~$K$~(see~\cite{P4}).
\bigskip

\begin{center}
\textbf{7. How to verify the  c.s.r. property~?}
\end{center}
\bigskip

Given a finite family of matrices~$\cA$, how to decide whether the semigroup $\cS_{\cA}$ is c.s.r~? Rather surprisingly, the answer is different for reducible and for irreducible
families. In the irreducible case, the c.s.r. property can be verified by an efficient polynomial time algorithm,
while in the reducible case, this problem is algorithmically undecidable even for a pair of nonnegative matrices. This means that there is no algorithm
that for any pair of nonnegative rational matrices $\cA = \{A_1, A_2\}$,  decides this property within finite time.
We start with the case of positively irreducible semigroups, when this problem is easily solvable.

\smallskip

\begin{center}
{\tt 7.1. Positively irreducible semigroups}
\end{center}
\smallskip

 To decide, whether  a given positively irreducible semigroup~$\cS$ is c.s.r. it suffices to take
a positively irreducible matrix $\bar A \in {\rm co}(\cS)$, which always exists, take its unique, up to normalization,
Perron-Frobenius eigenvector~$v$ and consider the smallest (by inclusion) common invariant linear subspace~$\tilde V$
that contains all the vectors~$v - Av\, , \, A \in \cS$. Then $\cS$ is c.s.r. if and only if  $v \notin \tilde V$.

If $\cS$ is generated by a finite family~$\cA = \{A_1, \ldots , A_m\}$, then one can take $\bar A = \frac1m \sum_{i=1}^m A_i$
and $\tilde V$ is the smallest common invariant subspace of~$\cA$ that contains  the vectors~$v - A_iv\, , \, i = 1, \ldots , m$. To find $\tilde V$, we denote by $\tilde V_1$ the linear span of those $m$ vectors and iteratively
construct the sequence~$\tilde V_1 \subset \tilde V_2 \subset \ldots $ as follows: $\tilde V_{k+1}$ is the linear span of $\tilde V_k$ and of $A_i\tilde V_k, \, i = 1, \ldots , m$. For the smallest $k\le d-1$ such that $\tilde V_k = \tilde V_{k+1}$, the subspace $\tilde V = \tilde V_k$ is desirable. The algorithm performs at most $d-1$ iterations,
each consists of verifying whether a given new vector $A_ib_s$ belongs to the current subspace or has to be added to the basis
($b_1, \ldots , b_r$ is a basis of a current subspace, $r\le d-1$). Each verification is by computing of the corresponding determinant,
and this has to be done for each $i = 1, \ldots , m$
and $s = 1, \ldots , r$. Thus, at most $d-1$ iterations, each involves at most  $m(d-1)$ determinant computations in dimensions at most~$d$. We see that the algorithm is polynomial.
\smallskip

The case of general irreducible semigroups (not necessarily nonnegative) is harder.
Nevertheless, at least two effective methods exist for this case.

\smallskip

\begin{center}
{\tt 7.2. Irreducible semigroups. The first method}
\end{center}
\smallskip

The algorithm for nonnegative semigroups can be extended directly to arbitrary irreducible semigroups by the Kronecker lifting.
To any $d\times d$-matrix $A$, one associates a linear operator $A^{\otimes \, 2} = A\otimes A$ on the space $\cM_d$
of symmetric $d\times d$-matrices defined as follows: $A^{\otimes \, 2} X \, = \, A^T X A\, , \  X \in \cM_d$.
This is nothing else but the tensor product of the matrix $A$ with itself. We have
$\rho(A^{\otimes  2}) \, = \, \rho^{\, 2}(A)$. Moreover, $A^{\otimes \, 2}$ leaves the cone~$\cK_d$
of positively semidefinite matrices invariant. Thus, the Kronecker lifting transfers any
family~$\cA$ of $d\times d$-matrices to the family~$\cA^{\otimes \, 2} = \{A^{\otimes \, 2} \ | \  A \in \cA\}$
of linear operators over the $(d^2+d)/2$-dimensional space $\cM_d$
with a common invariant cone~$\cK_d$, and the same is for semigroups $\cS \, \to \, \cS^{\, \otimes \, 2}$.
By increasing the dimension, we obtain a semigroup with a common invariant cone and with the same spectral properties.
In particular, $\cS$ is c.s.r. if and only if $\cS^{\otimes \, 2}$ is. Moreover, $\cS$ is irreducible
if and only if $\cS^{\times \, 2}$ is positively irreducible with respect to the cone~$\cK_d$, i.e.,
none faces of~$\cK_d$ is a common invariant subspace for~$\cS^{\, \otimes \, 2}$. Therefore, we can verify the c.s.r. property
of $\cS^{\otimes \, 2}$ by the same algorithm as for nonnegative matrices, using an analog
of Proposition~\ref{p50} for the cone $\cK_d$ instead of $\re^d_+$ (Remark~\ref{r30}). This is the idea of the
first method of deciding the c.s.r. for irreducible family.
\smallskip

{\em The algorithm}. Given an irreducible semigroup~$\cS$, we need to decide
whether it has the constant spectral radius or not.
We take a positively irreducible operator $\widehat A \in {\rm co}(\cS^{\otimes \, 2})$, take its
Perron-Frobenius eigenvector~$Y \in \cK_d$ and consider the smallest (by inclusion) common invariant linear
subspace~$\tilde V \subset \cM_d$
that contains all the matrices~$Y - A^TYA\, , \, A \in \cS$. Then $\cS$ is c.s.r. if and only if  $Y \notin \tilde V$.
{\em End}.
\smallskip

If $\cS$ is generated by a finite family~$\cA = \{A_i\}_{i=1}^m$, then one can take
$\widehat A \, = \, \frac1m \sum_{i=1}^m A_i^{\otimes \, 2}$
and $\tilde V \subset \cM_d$ is the smallest common invariant subspace of~$\cA^{\otimes \, 2}$
that contains  the matrices ${Y - A_i^TYA_i\, , \, i = 1, \ldots , m}$. The construction of the
subspace $\tilde V$ is the same as for semigroups of nonnegative matrices, the complexity estimate
is also the same, with replacement of the dimension from $d$ to $(d^2 + d)/2$. Thus,
if all the matrices $A_i$ are rational, then the subspace~$\tilde V$ is found and the relation
$Y \notin \tilde V$ is verified within polynomial time in the input. This proves the following theorem:
\begin{theorem}\label{th60}
The c.s.r. recognition problem is solved by a polynomial time algorithm for an irreducible
family of rational matrices, as well as for a positively irreducible family of nonnegative rational matrices.
\end{theorem}

\begin{remark}\label{r40}
{\em In fact, the common affine invariant subspace $V$ for all operators~$A^{\otimes \, 2}, \, A \in  \cS$
is closely related to walks on the Euclidean sphere from Section~5.  By Theorem~\ref{th20},
an irreducible semigroup~$\cS$ is c.s.r. precisely when there is a point $x \in \re^d$
and a matrix $H \in \cK_d$ such that $x^THx = 1$ and $x^T(A^THA)x = 1$ for all $A \in \cS$.
Indeed, after linear change of coordinates taking the ellipsoid $\{y \in \re^d \ | \ y^THy = 1\}$
to the unit Euclidean ball, these assertions become $\|x\|_2 = 1, \,  \|Ax\|_2 = 1, \, A \in \cS$.
Hence, the affine hull~$V$ of all the matrices  $A^THA, \, A \in \cS$, is invariant with respect to
the operators of the semigroup~$\cS^{\, \otimes \, 2}$
and is contained in the affine  hyperplane $V_x = \{M \in \cM_d \ | \  x^TMx = 1\}$ of $\cM_d$.
Thus, $V$ is defined by the orbit of the point $H \in \cM_d$, which is the matrix
of the ellipsoid, whose surface contains all the walk points~$Ax, \,  A \in \cS$.
Every point~$H$ from the intersection of $V$ with the interior of~$\cK_d$ is the quadratic form
of such an ellipsoid, and all of them can be found by mere evaluating of the invariant subspace~$V$.
The point $x\in \re^d$, in turn, defines the desired affine hyperplane~$V_x \subset \cM_d$ which contains~$V$.
In practice, however, this is not clear how to find this point $x$ algorithmically. The problem is apparently
reduced to a system of quadratic equations of $d$ variables. That is
 why, in the algorithm above we construct the invariant subspace~$V$ in a different way.
}
\end{remark}

\smallskip

\begin{center}
{\tt 7.3. Irreducible semigroups. The second method}
\end{center}
\smallskip

The method is based on the fact that the c.s.r. property of a semigroup $\cS$
generated by an irreducible finite family~$\cA = \{A_i\}_{i=1}^m$ is equivalent to
the assertion $\rho_2 (\cA) = \rho_4(\cA) = 1$, where $\rho_p$ is the $p$-radius (Definition~\ref{d8}).
The proof is actually the same as for the assertion $\rho_1(\cA) = \rho_2(\cA) = 1$
from the proof of Theorem~\ref{th20}. The values $\rho_2 (\cA)$ and $\rho_4 (\cA)$ are both
efficiently computed: $\rho_2(\cA)$ is the square root of the spectral radius of the operator
$\frac1m \sum_{i=1}^m A_i^{\otimes 2}$, and  $\rho_4(\cA)$ is the power $1/4$ of the spectral radius of the operator
$\frac1m \sum_{i=1}^m A_i^{\otimes 4}$, where $A_i^{\otimes 4} = (A_i^{\otimes 2})^{\otimes 2}$
is the operator in the space of dimension~${d \choose 4}$.
\medskip

\smallskip
\newpage

\begin{center}
{\tt 7.4. General semigroups. Undecidabilty}
\end{center}
\smallskip

For reducible semigroups, however, the c.s.r. recognition problem harder.
This is explained by Theorem~\ref{th5}: to verify the c.s.r. property one needs to check that
the joint spectral radii of the blocks $\cS^{(i)}$ in factorization~(\ref{blocks}) do not exceed one.
However, the problem of verifying the inequality $\rho(\cA) \le 1$ is algorithmically undecidable even for a pair
of nonnegative rational matrices~$\cA$ (see~\cite{BT, TB}).
\begin{theorem}\label{th70}
The c.s.r. recognition problem for a pair of nonnegative  rational matrices is algorithmically undecidable.
\end{theorem}
{\tt Proof}. To an arbitrary pair~$\cA = \{A_1, A_2\}$ of $d \times d$ nonnegative matrices we associate
a pair $\cA' = \{A_1', A_2'\}$ of $(d+1)\times (d+1)$ matrices as follows:
$A_i'$ has two diagonal blocks, the first block is $\{1\}$ (one-dimensional) the second block is
$A_i$ ($d$-dimensional), all other entries of $A_i'$ are zeros. By Theorem~\ref{th30}, the semigroup
generated by the pair $\cA'$  has constant spectral radius if and only if $\rho(\cA) \le 1$. Since the latter assertion is
algorithmically undecidable, the former also is.

{\hfill $\Box$}

\smallskip

\medskip

\begin{center}
\textbf{8. A classification of c.s.r. semigroups}
\end{center}
\bigskip

Apart from the case of nonnegative matrices (Section~6), the problem of complete classification
of c.s.r. semigroups is not solved. In the next subsection, we describe several classes of c.s.r semigroups
 and leave an open question about the completeness of that list. Then we prove Lemma~\ref{th10} on the structure of s.c.r.
 semigroups, using which we classify them in low dimensions: $d=2, 3$.
 By Theorem~\ref{th5} it suffices to consider the irreducible case only. In what follows in this section, all semigroups are assumed to be irreducible.  The main parameter of the classification
 is $r(\cS)$, the smallest rank of matrices from~$\cS$. The case $r = d$ is simple: $\cS$ consists of orthogonal matrices
 (Theorem~\ref{th10}).
 The case $r=d-1$ is also classified, in all other cases we do not have a complete answer for $d \ge 4$.
 This issue is discussed in Subsection~8.5.

\begin{center}
{\tt 8.1. The list  of irreducible c.s.r. semigroups}
\end{center}
\bigskip

We spot the following classes of c.s.r. semigroups:
\medskip

\noindent \textbf{1.} {\em An arbitrary  subgroup of $O(d)$.}
\smallskip

\noindent \textbf{2.} {\em An arbitrary semigroup of $(k, n)$-matrices}.
\smallskip

For  an irreducible c.s.r. semigroup~$\cS_k$
of $k\times k$ matrices, $k \ge 1$, and for an arbitrary $n \ge 2$, a $(k,n)$-matrix~$A$ is defined as a matrix
of size $kn$ that consists of $n^2$ square blocks of size~$k$; each of $n$ block rows of~$A$ contains a unique nonzero
block from~$\cS_k$, all other blocks in that row are zero.  Thus, any $(k, n)$-matrix contains exactly
$n$ nonzero blocks of size $k$, each of them equals to some matrix from~$\cS_k$. Clearly, any product of
$(k,n)$-matrices (corresponding to the same semigroup~$\cS_k$) is again a $(k, n)$-matrix.
 Hence, the norms of products are bounded uniformly from above and from below.
Therefore, every set of $(k,n)$-matrices generate a c.s.r. semigroup of matrices of size $d = kn$.

Thus, for every irreducible c.s.r. semigroup of $k\times k$-matrices and for any $n \ge 2$,
the transfer to $(k, n)$-matrices produces c.s.r. semigroups of dimension~$kn$.

\begin{ex}\label{ex30}
{\em An $(1, n)$-matrix is a matrix of size~$n$ that has a unique nonzero element $\pm 1$ in each row, all other elements
are zeros. If we denote by $\cI_n = \{x \in \re^n \ | \ x_i = \pm 1\, , \, i = 1, \ldots , n\}$
the set of vertices of an $n$-dimensional cube, then any $(1, n)$-matrix~$A$ satisfies $A\, \cI_n \, \subset \cI_n$.
This is a characteristic property of $(1, n)$-matrices.
}
\end{ex}

\smallskip

\noindent \textbf{3.} {\em An arbitrary semigroup of $(\bk, \bn)$-torsion matrices}.

Let $s \ge 2, \bk = (k_1, \ldots , k_s) \in \z^{s}_+, \bn = (n_1, \ldots , n_s)\in \z^s_+$,
where $k_i, n_i \ge 1$ for all $i$ and $n_j \ge 2$ for at least one~$j$. We define two matrices $B, C$
of size $d = \sum_{i=1}^s k_in_i$ as follows. The matrix $C$ is block-diagonal, it contains
$s$ square diagonal blocks~$C_i$ of sizes $k_in_i, i = 1, \ldots , s,$ all other entries of $C$ are zero.
The $i$th bock~$C_i$ is also a block matrix, it is composed of $n_i^2$ square blocks of size~$k_i$.
Exactly one of $n_i$ block columns of $C_i$ is nonzero, denote its number by $h_i$.
All the $n_i$ blocks $Q_{i, h_i, 1}, \ldots , Q_{i, h_i, n_i}$ of this column are
some orthogonal $k_i\times k_i$-matrices. Let $L_{ij}\subset \re^d$ be the $k_i$-dimensional subspace
corresponding to the $j$th block column of $C_i$. Thus, $\re^d \, = \, \oplus \sum_{i, j}L_{i,j}$
and, respectively, $x = \sum_{i, j}x_{i,j}\, , \ x_{ij} \in L_{ij}$.
Consider the following subspace of $\re^d$:
$$
L\  = \ \Bigl\{\, x \in \re^d \ | \  Q_{i, h_i, 1}^{-1}x_{i1} = \cdots = Q_{i, h_i, n_i}^{-1}x_{in_i} \, , \,
i = 1, \ldots , s\, \Bigr\}
$$
This subspace contains the range of $C$. Finally, let $B$ be an arbitrary orthogonal operator on~$L$.
Then $A = BC$ is a $(\bk, \bn)$-torsion matrix.

Denote the following  subset of the unit sphere in $\re^d$:
$$
\cI_{\, \bk, \bn} \ = \ \Bigl\{\, x \in \re^d \ | \ \|x_{i1}\|_2 = \cdots = \|x_{in_i}\|_2 , \, \,
i = 1, \ldots , s\, , \, \|x\|_2 = 1\Bigr\}\, .
$$
For any $(\bk, \bn)$-torsion matrix~$A$ (the vectors~$\bk, \bn$ are fixed, the indices $h_i \in \{1, \ldots , n_i\}$ and the orthogonal matrices $Q_{i, h_i, j}$ are arbitrary), we have $A\cI_{\, \bk, \bn} \subset \cI_{\, \bk, \bn}$.
Hence, all products of $(\bk, \bn)$-torsion matrices define walks on the Euclidean sphere, and hence, each set of
those matrices defines a c.s.r. semigroup.
\smallskip

\smallskip

\noindent \textbf{4.} {\em An arbitrary semigroup generated by tensor products $A_1\otimes A_2$},
where $A_i \in \cS_i$, and $\cS_i$ is an irreducible  c.s.r. semigroup  of matrices of size~$d_i \ge 2\, , \, i = 1, 2$.
\smallskip

Thus, having a pair of irreducible c.s.r. semigroups $\cS_1, \cS_2$ of dimensions $d_1$ and $d_2$ respectively, we take an arbitrary set of tensor
 products of their matrices. This set generates a c.s.r. semigroup in dimension~$d_1d_2$.
\smallskip

\noindent \textbf{5.} {\em The set of transpose matrices to a c.s.r. semigroup}.
\smallskip

In particular, the semigroups of $(k,n)$-matrices from the class~2 defined by columns instead of rows
(there is a unique nonzero block in each block column) is also c.s.r.

Of course, the linear similarity transform (a transfer to another basis) respects the c.s.r. property. We do not take it into account and
consider all the classes 1 -- 5 up to a linear similarity.    

\medskip

\noindent \textbf{Problem 1}. {\em Is it true that every irreducible c.s.r. semigroup is obtained from
a subgroup~of~$O(d)$ (item 1) by the four procedures from items 2 -- 5 and by linear similarities~?}
\smallskip

In Subsections 8.3 and 8.4,  we shall see that the answer is affirmative for $d=2,3$ and for some other special cases. To attack  the low-dimensional cases, we prove the following lemma which gives the main tool for c.s.r. classification.
\bigskip

\begin{center}
{\tt 8.2. The key lemma}
\end{center}
\bigskip

We write $r(\cS)$ for the {\em rank of the semigroup}~$\cS$,
i.e., the minimal rank of matrices from~$\cS$. Since $\cS$ is of constant spectral radius, it follows that
$ r \ge 1$. Theorem~\ref{th10} implies that if $r = d$ for an irreducible semigroup~$S$, then it has a constant spectral radius if and only if~$\cS$ consists of orthogonal
matrices (in a suitable basis). Let us note that if $\cS$ is irreducible and c.s.r., then
its rank coincides with the rank of its closure. Indeed, every matrix from $\cS$ has eigenvalues either zeros or
of modulo~one~\cite{OR}, hence taking a limit of such matrices does not reduce the rank.
Without loss of generality we assume below that~$\cS$ is closed.
\begin{lemma}\label{l10}
For any irreducible c.s.r. semigroup~$\cS$, there is a subspace $\, L \subset \re^d,
\, {\rm dim}\, L\, = \, r$, and a basis in~$\re^d$ such that the orthogonal projector~$P$ onto $L$ belongs to~$\cS$
and for every $A \in \cS$, the operator $PA\,|_{L}$ is orthogonal.
 \end{lemma}
 The proof is in Appendix.

\begin{cor}\label{cor20}
 For every irreducible c.s.r. semigroup~$\cS$ there is basis in~$\re^d$
 in which every matrix from~$\cS$ has its upper principal $r\times r$ submatrix
 (in the first $r$ rows and first $r$ columns) orthogonal.
 \end{cor}

Lemma~\ref{l10} has the following geometrical meaning. If $\, \fB_r$ is an $r$-dimensional Euclidean ball in~$L$,
then $PA(\fB_r) = \fB_r$ for every $A \in \cS$. In other words, every operator from $\cS$ maps the ball
$\fB_r$ to an $r$-dimensional ellipsoid, which is a cross-section of the right circular cylinder $\fC_r =
\{x + y\  | \  x \in \fB_r\, , \, y \perp L\}$ by some $r$-dimensional subspace of~$\re^d$.

\bigskip

\begin{center}
{\tt 8.3. Low dimensions:} ${d =2}$.
\end{center}
\bigskip

  If $r=2$, then by Theorem~\ref{th10} all matrices of $\cS$ are orthogonal. If $r=1$,
then by Lemma~\ref{l10} there is an orthogonal projector $P$ onto some one-dimensional subspace $L$
such that for every $A \in \cS$ the operator $PA$ is orthogonal on~$L$. Denote by $\fB = [b_1, b_2]$ the segment of the line
$L$ of length~$2$ centered at the origin. We have $PA \fB = \fB$ for each $A \in \cS$, hence
$A\fB$ is a segment centered at the origin with ends on the lines $J_1$ and $J_2$ that are orthogonal
to $\fB$ and pass through the points $b_1$ and $b_2$ respectively.  Thus, all operators $A \in \cS$ map $\fB$
to such segments. Let us show that the set $\{A\fB \ | \ A \in \cS\}$ actually consists of two segments, including
$\fB$. By Lemma~\ref{l10}, $\, P \in \cS$, hence the segment~$\fB$ belongs to this set. If all segments of this set
coincide with $\fB$, then $\cS$ is reducible. Hence, there is a segment $[c_1, c_2] \ne [b_1, b_2]$ in this set.
Assume there is a third  segment $[d_1, d_2]$, where $c_i, d_i \in L_i\, , \ i = 1,2$. If for all $A \in \cS$
the points $Ac_1$ and $Ad_1$ lie on the same line~$J_i$ then $(AJ_1)\perp L$ for all $A \in \cS$, hence $\cS$ has a common invariant subspace parallel to~$J_1$. Therefore, for some $A \in \cS$ the points $Ac_1$ and $Ad_1$ lie on two different  lines~$J_i$. However, in this case the point $Ab_1$ does not lie on these lines (since $Ab_1, Ac_1$ and $Ad_1$ are co-linear), which is impossible. Thus, the set of all images of $\fB$
consists of exactly two segments. Every operator $A \in \cS$ maps each of these segments either to itself or to the other one. After a suitable linear transform, we assume that these segments are two diagonals of a unit square~$\fQ$.
Then every operator from~$\cS$ is either a projection to one of diagonal of~$\fQ$ parallel to one of its sides,
or an isometry of~$\fQ$, or composition of several such operators. Thus, we have proved
\begin{prop}\label{p10}
All irreducible c.s.r. semigroups $\cS$ of $2\times 2$-matrices are classified as follows:

 if $r=2$, then $\cS$ consists of orthogonal matrices;

if $r=1$, then $\cS$ consists of compositions of the following operators:

1) projections on a diagonal of $\fQ$ parallel to its side ($4$ operators);

2) orthogonal transforms of~$\fQ$ (two rotations on~$\pi/2$, four axial symmetries, and the central symmetry).

\end{prop}

If~$\fQ = \{x \in \re^2 \ | \ \max \{|x_1|, |x_2|\} \le 1\} $,  then the matrices
from classes 1) and 2) in Proposition~\ref{p10} are actually $(1,2)$-matrices (Example~\ref{ex30}).
Thus, in  the case~$d  = 2, r = 1$, every irreducible c.s.r. semigroup
consists of matrices with exactly one nonzero element $\pm 1$ in each row.

\begin{cor}\label{c20}
If an irreducible semigroup of $2\times 2$-matrices of constant spectral radius
contains at least one degenerate matrix, then it is finite and consists of
$(1,2)$-matrices.
\end{cor}
\medskip

\begin{center}
{\tt 8.4. Low dimensions:} ${d =3}$.
\end{center}
\bigskip

\noindent  If $r=3$, then all matrices of $\cS$ are orthogonal. Consider the other cases.
\smallskip

\textbf{The case} $\mathbf{r=1}$. By Lemma~\ref{l10} there is an orthogonal projector $P$ onto some one-dimensional subspace $L \subset \re^3$ such that for every $A \in \cS$ the operator $PA$ is orthogonal on~$L$. Let $\fB = [b_1, b_2]$ be the segment of $L$ of length~$2$ centered at the origin. For each $A \in \cS$, the segment
$[x_1, x_2] = A\fB$ is centered at the origin and $x_i \in J_i\, , \ i = 1,2$ where $J_1$ and $J_2$
are two-dimensional planes  orthogonal
to $\fB$ and passing  through the points $b_1$ and $b_2$ respectively.
From the irreducibility of~$\cS$ it follows that for some $C \in \cS$ the corresponding segment $[c_1, c_2] = C\fB$
does not coincide with $[b_1, b_2]$. Moreover, there is an operator $D \in \cS$ that sends
the points $b_1$ and $c_1$ to different planes~$L_i$. Otherwise the linear span of the vectors $\{A(c_1 - b_1) \ | A \in \cS\}$ is a common nontrivial invariant subspace for~$\cS$.
Take now arbitrary $A \in \cS$ and consider the segment
$[x_1, x_2] = A\fB$. The point $Dx_1$ belongs to either $J_1$ or $J_2$. The set $\{ x \in J_1 \ | \ Dx \in J_1\}$
is either empty, or a line in~$J_1$, of the whole~$J_1$. Since it
contains precisely one of the two points $b_1$ and $c_1$, it is neither empty nor the whole~$J_1$, so it is a line.
The set $\{ x \in J_1 \ | \ Dx \in J_2\}$ is a parallel line. Thus, $x_1$ belongs to one of these parallel lines.
Since this is true for every $A\in \cS$, we see that the ends of all segments $A\fB\, , \ A \in \cS$ on the plane~$J_1$ lie on this pair of parallel lines. By the irreducibility of~$\cS$, there is another pair of parallel lines, which is not parallel
to the first pair, that  also contains all these ends of the segments $A\fB\, , \ A \in \cS$. Whence, all these ends
are located at common points of these four lines, i.e., at
four vertices of a parallelogram. Therefore, all the segments $A\fB\, , \ A \in \cS$ are diagonals of the parallelepiped~$\fP$ centered at the origin, whose face coincide with that parallelogram. After a suitable linear transform, it may be assumed that~$\fP$ is a unit cube, whose vertices have coordinates $\pm 1$. If a matrix $A \in \cS$
maps each diagonal of~$\fP$ to a diagonal, then $A$ is a $(1,3)$-matrix (Example~\ref{ex30}).
If $A$ respects the set of three diagonals, then, in the basis of these diagonals, $A$ is a transpose
to a $(1,3)$-matrix. Thus, we have proved
\begin{prop}\label{p20}
Up to a linear similarity, every irreducible c.s.r. semigroup of $3\times 3$-matrices for which $r=1$,
either consists of $(1,3)$-matrices, or consists of transposes to $(1,3)$-matrices.
\end{prop}
Thus, in case $d=3, r=1$, every irreducible c.s.r. semigroup either consists of matrices with a unique
nonzero element $\pm 1$ in each row, or consists of matrices with a unique nonzero element $\pm 1$
in each row.
\bigskip

\textbf{The case} $\mathbf{r=2}$. By Lemma~\ref{l10}  there is an orthogonal projector $P$ onto a two-dimensional subspace $L \subset \re^3$ such that for every $A \in \cS$ the operator $PA$ is orthogonal on~$L$. Let $\fB$ be the
Euclidean unit ball in $L$. For each $A \in \cS$, the set $A\fB$ is a cross-section of a
right circular cylinder $\fC =
\{x + y\  | \  x \in \fB\, , \, y \perp L\}$ by some two-dimensional subspace of~$\re^3$.
Assume there are operators $A_1, A_2 \in \cS$ such that the three sets $A_1\fB, A_2\fB$, and $A\fB$ are all
different. Let a line orthogonal to the plane~$\fB$ on the surface of the cylinder meets these cross-sections
at points~$a_1, a_2$, and~$a$ respectively. From the irreducibility of~$\cS$ it follows that there is an operator~$D \in \cS$
such that the line passing through $Da_1$ and $Da_2$ is not orthogonal to~$L$. Hence, this line intersects the surface of
the cylinder~$\fP$ at most at two points, which is impossible, because the three points $Da_1, Da_2, Da$ belongs to that
intersection. Thus, the set $\{A\fB \ | \ A \in \cS\}$ consists of two cross-sections of the cylinder, one of which is~$\fB$.
Up to an affine similarity, it may be assumed that these two cross-sections are equal concentric Euclidean discs orthogonal to each other.
The semigroup~$\cS$ consists of operators that map each of these discs to itself or to another disc.

\begin{prop}\label{p30}
Up to a linear similarity, every irreducible c.s.r. semigroup of $3\times 3$-matrices for which $r=2$,
consists of operators that map each of the two equal concentric orthogonal  Euclidean discs to itself or to the other disc.
\end{prop}
Such semigroup may contain the following operators and their compositions:

 1) a composition of the projection of the first disc onto the second one with an orthogonal transform of the second disc;

2) orthogonal transforms respecting the discs:   rotation on~$\pi/2$ and on~$\pi$ around  the common line of the discs,
symmetries with respect to the planes of the discs and with respect to their bisector planes, central symmetry with respect to the origin.
\smallskip
\enlargethispage{\baselineskip}
If the planes of the disc have the equations $L_{\pm} = \{x \in \re^3 \ | \  x_1 \, = \, \pm \, x_2\}$,
then the operators of the class 1) are given by $(\bk, \bn)$-torsion matrices with $\bk = (1, 1), \, \bn = (2, 1)$
(item 3 in the list of Subsection~8.1).
Each of them is the product $A = BC$, where $C$ consists of two diagonal blocks: $C_1$ is a $2\times 2$ matrix
with one nonzero column with entries $\pm 1$; $C_2$ is $\pm 1$. The matrix $B$ defines an orthogonal transform
of the two-dimensional plane~$L_{+}$ or $L_{-}$ (the plane containing the disc).

The operators of class 2) are given by $(1,3)$-matrices (Example~\ref{ex30}).

\bigskip

\medskip

\begin{center}
{8.5.  Arbitrary dimension~$d$}
\end{center}
\bigskip

Already for $d=4$ we are not able to classify all irreducible c.s.r. semigroups. Using items 1 -- 5 in Subsection~8.1
we obtain the following c.s.r. families of $4\times 4$-matrices:
\smallskip

an arbitrary subset of the group $O(4)$;
\smallskip

an arbitrary set of $(1,4)$-matrices; a transpose to that set;
\smallskip

a set of $(2, 2)$-matrices made of an arbitrary  c.s.r family of matrices of size~$2$,
a transpose to that set;
\smallskip

an arbitrary set of tensor products of two c.s.r. families of $2\times 2$ matrices (Subsection~8.3);
\smallskip

an arbitrary set of $(\bk, \bn)$-torsion matrices (item 3, Subsection~8.1)
 either with $\bk = (1, 1), \, \bn = (3, 1)$, or with $\bk = (1, 1), \, \bn = (2, 2)$, or with
$\bk = (1, 2), \, \bn = (2, 1)$, or with $\bk = (1, 1, 1), \, \bn = (2, 1,1)$.
\medskip

Of course, all families obtained by those ones by a change of basis are also c.s.r. We do not know
whether this list is complete.
   \medskip

For arbitrary $d$, the construction of c.s.r. semigroups is realized by items 1 -- 5 of Subsection~8.1 in the same way,
although the number of cases grows dramatically  with the dimension. The problem of completeness of this classification
is solved by now only for two values of~$r$. If $r = d$, then, apparently, $\cS$ is a subgroup of~$O(d)$.
If $r=d-1$, then a complete analogue of Proposition~\ref{p30} takes place (with essentially the same proof):
every irreducible c.s.r. semigroup of $d\times d$-matrices for which $r=d-1$,
consists of operators that map each of the two equal concentric orthogonal  Euclidean $(d-1)$-dimensional balls to itself or to the other ball. The class 1) of Proposition~\ref{p30} is described now
by $(\bk, \bn)$-torsion matrices with $\bk = (1, \ldots , 1), \, \bn = (2, 1, \ldots , 1)$
(the vectors are both $(d-1)$-demensional); the class 1) is by $(1, d)$-matrices.

\bigskip

\newpage 

\begin{center}

\textbf{9. Applications}

\end{center}

\bigskip

\smallskip
\begin{center}
{\tt 9.1. Finite matrix semigroups}
\end{center}
\smallskip

Finite matrix semigroups have been studied in the literature, see~\cite{CCJ, Jacob, McNZ, MS} and bibliography in those works. In fact, they are
closely related to c.s.r. semigroups.
\begin{prop}\label{p90}
If a finite matrix semigroup~$\cS$ does not contain zero matrix, then~$\cS$ is c.s.r.
\end{prop}
{\tt Proof.} If, to the contrary, $\rho(A) \ne 1$ for some $A \in \cS$, then
all numbers $\rho(A^k) = \rho^k(A)$ are different, whenever~$\rho(A) \ne 0$, and therefore, $\cS$
contains an infinite subset~$\{A^k, \ k \in \n\}$. If $\rho(A) = 0$, then $A$ is nilpotent,
and hence $A^d = 0$, which is impossible, since $A^d \in \cS$.

 {\hfill $\Box$}
\smallskip

Invoking Theorem~\ref{th10}, we obtain
\begin{cor}\label{c60}
For every irreducible finite matrix semigroup~$\cS$, there is a basis in~$\re^d$ in which all
nonsingular elements of~$\cS$ are orthogonal. Every irreducible finite semigroup of
nonsingular matrices is, in a suitable basis, a subgroup of~$O(d)$.
\end{cor}
{\tt Proof.} The first part follows from Theorem~\ref{th10}. In the second part, it remains to prove that
a finite semigroup~$\cS$ of nonsingular matrices is a group, i.e., contains an inverse to each~$A \in \cS$.
Indeed, since the sequence $\{A^k\}_{k \in \n}\subset \cS$ is finite,
it follows that $A^{k} = A^{k+n}$ for some $k, n \in \n$. Hence, $A^n = I$ and $A^{n-1} = A^{-1}$.
Since $A^{n-1} \in \cS$, the proof is completed.

 {\hfill $\Box$}
\smallskip

Thus, all finite irreducible semigroups of nonsingular matrices are finite subgroups of~$O(d)$.
Note that every nonsingular element of a finite matrix semigroup is evidently similar to
an orthogonal matrix, otherwise all its powers are different. What is nontrivial in Corollary~\ref{c60}
is that all nonsingular matrices become orthogonal simultaneously, in one basis.

Finite semigroups not containing zero matrices form a subclass of c.s.r. semigroups.
For integer matrices, i.e., for matrices with integer entries, these classes coincide.
\begin{prop}\label{p110}
An irreducible c.s.r. semigroup of integer matrices is finite.
The same holds for positively irreducible semigroups of nonnegative integer matrices.
\end{prop}
{\tt Proof.} All norms of matrices in an irreducible c.s.r. semigroup are bounded by some constant, hence
the number of integer matrices is finite.

 {\hfill $\Box$}
\smallskip

The problem of classification of finite semigroups is still open, even for integer matrices.
The results of Section~8 give such a classification in low dimensions ($d=2,3$) and many examples
of finite semigroups in higher dimensions. Starting with a finite subgroup of~$O(d)$ and applying
procedures 2 -- 5 from the list in Subsection~8.1 we obtain finite semigroups. As in Problem~1, one may rise the
question about the completeness of this classification.

The algorithmic recognition of finiteness for the semigroup generated by a given family of matrices
$\cA = \{A_1, \ldots , A_m\}$ is also a challenging problem. Algorithms presented in~\cite{MS, JPB}
solve it for nonnegative integer matrices (the algorithm from~\cite{JPB} is polynomial). We make the next step and solve the problem  for arbitrary irreducible family~$\cA$ of integer matrices, provided it is not mortal, i.e., none of products
of those matrices is zero. Note that the mortality recognition problem is known to be
algorithmically undecidable~\cite{TB}.
\begin{theorem}\label{th110}
There is a polynomial-time algorithm that for every irreducible non-mortal family
of integer matrices $\cA = \{A_1, \ldots , A_m\}$ decides whether it generates a finite semigroup.
\end{theorem}
{\tt Proof}. In view of Propositions~\ref{p90} and~\ref{p110}, the finiteness of
$\cS_{\cA}$ is equivalent to its c.s.r. property, which can be decided by the algorithm from Subsection~7.
(Theorem~\ref{th60}).

 {\hfill $\Box$}
\smallskip

There are many examples of finite semigroups of integer matrices.
We consider one class obtained by iterating the construction of
  $(k, n)$-matrices (item 3, Subsection~8.1).
\begin{ex}\label{ex50}
{\em Take a finite sequence
$n_1, \ldots , n_q$, each $n_i$ is natural and is not one. Consider an arbitrary set
of $(1, n_1)$-matrices (Example~\ref{ex30}) and transpose them. We get a set of
matrices of size~$n_1$, each of them has one nonzero element $\pm 1$ in every column.
Then, from this set, construct an arbitrary set of $(n_1, n_2)$-matrices (item 2 in the list in Subsection~8.1)
and transpose them. Then construct a set of $(n_1n_2\, , n_3)$-matrices, transpose them, and so on.
As a result, we obtain a c.s.r. family~$\cA$ of matrices of size $d = n_1\ldots n_q$.
All their entries are zeros and plus/minus ones, all their products are uniformly bounded, hence
the semigroup~$\cS_{\cA}$ is finite. Take now an arbitrary integer matrix~$C$
such that $|{\rm det}\, C| = 1$. Then the family $\{C^{-1}AC \ | \  A \in \cA\}$ generates
a finite semigroup of integer matrices.
}
\end{ex}

\smallskip
\begin{center}
{\tt 9.2. Linear switching systems}
\end{center}
\smallskip

Given a compact family $\cA$ of $d\times d$ matrices, a linear switching system (LSS) is the following linear differential
equation on the vector-valued function $x: [0, +\infty) \to \re^d$:
\begin{equation}\label{cont}
\left\{
\begin{array}{l}
\dot x (t) \ = \ A(t)\, x(t)\,  ;\\
x(0) \, = \, x_0\, ,
\end{array}
\right.
\end{equation}
where $A(\cdot): \, [0, +\infty) \to \cA$ is a measurable function called the {\em switching law}.
The solution~$x(\cdot)$ is a {\em trajectory} of the system corresponding to the  switching law~$A(\cdot)$ and to the
initial condition $x(0) = x_0$. There is an extensive literature on LSS, in particular,
exploring the asymptotic growth of the trajectories~(see~\cite{Lib, LA, MP} and references therein).
 The {\em exponent of growth} of the switching law $A(\cdot)$ is
$$
\sigma (A(\cdot)) \ = \ \sup_{x_0 \in \re^d}\,  \Bigl( \, \limsup_{t \to \infty} \, \frac1t \,  \log \, \|x(t)\|\, \Bigr) \, ,
$$
where~$x(t)$ is the trajectory with the initial condition~$x(0) = x_0$. This is the fastest growth of all trajectories
 corresponding to a given switching law. We call a system {\em uniform} if $\sigma(A(\cdot))$ is the same for all
 switching laws $A(\cdot)$. Our goal is to characterize uniform systems.

 Observe that making an $\alpha$-shift of the system: $\, \cA \, \mapsto \, \cA + \alpha I = \{A + \alpha I \ |\ A \in \cA\}$,
 we replace every trajectory $x(t)$ by $e^{\alpha t} x(t)$, and, therefore, add $\alpha$ to every exponent
 of growth  $\sigma(A(\cdot))$. Hence, up to a suitable shift, it may be assumed that $\sigma (A(\cdot)) = 0$,
 and it suffices  to characterize only such uniform systems.
 \begin{prop}\label{p70}
 An LSS is uniform with $\sigma (A(\cdot)) = 0$ for every switching law, if and only if the family
 $\cM = \{e^{t A} \ | \ A \in \cA \, , \, t \in \re_+\}$ is c.s.r.
 \end{prop}
{\tt Proof.} If  $\cM$ is not c.s.r., then there are numbers $t_i > 0$ and matrices $A_i \in \cA, i = 1, \ldots , n$,
such that $\lambda = \rho(\prod_{i=1}^n e^{t_iA_i}) \ne 1$. Denote $T = \sum_{i=1}^n t_i$.
For the piecewise-constant switching law $A(t)$
taking values~$A_n, \ldots , A_1$ on successive segments of lengths $t_n, \ldots , t_1$,
we have $x(T) =  e^{t_1A_1}\cdots e^{t_nA_n}x(0) $. Therefore, for the periodization of this switching law
with period $T$, we have $\sigma (A(\cdot)) = T^{-1}\log \, \lambda \ne 0$, hence the LSS is not uniform.
For the converse, in view of Theorem~\ref{th5}, it suffices to consider an irreducible family~$\cA$, in which case $\cM$ is also irreducible. If $\cM$ is c.s.r., then the norms of all products of its matrices
is between two positive constants $C_1$ and $C_2$. Hence, for any piecewise-constant switching law, we have
$$
C_1 \ \le \ \max_{\|x_0\| = 1}\, \|x(t)\| \ \le \  C_2\, .
$$
 By continuity, this holds for every switching law~$A(\cdot)$,
and hence $\sigma (A(\cdot)) = 0$.

 {\hfill $\Box$}
\smallskip

By Theorem~\ref{th5}, it suffices to characterize only irreducible uniform LSS
(corresponding to irreducible families~$\cA$). Note that all matrices of~$\cM$ are invertible,
hence the c.s.r. property can be characterized by applying Theorem~\ref{th10}.
\begin{theorem}\label{th80}
 An irreducible  LSS is uniform with $\sigma = 0$ if and only if there is a basis in~$\re^d$
 in which all matrices of~$\cA$ are antisymmetric ($A^T = -A$). In this case, every trajectory~$x(t), \, t \in [0, +\infty)$,
 lies on a Euclidean sphere.
\end{theorem}
{\tt Proof}. If the LSS is uniform, then the family
of nonsingular matrices~$\cM$ is c.s.r., and by Theorem~\ref{th10}, all
the matrices $e^{tA}$ are orthogonal in a suitable basis. This means $A^T = -A$ for all~$A \in \cA$.
Conversely, if all matrices of~$\cA$ are antisymmetric, then $\bigl( \dot x(t)\, , \, x(t)\bigr) \, = \,
\bigl( A(t) x(t)\, , \, x(t)\bigr) = 0$, and consequently
$\|x(t)\|_2 \equiv {\rm const}$ along every trajectory.

 {\hfill $\Box$}
\smallskip

An important class of LSS are {\em positive systems}, for which all trajectories are nonnegative
($x(t) \ge 0\, , \, t \in \re_+$), provided $x(0) \ge 0$.
An LSS is positive precisely when each matrix $A \in \cA$ is Metzler, i.e. all off-diagonal
elements of~$A$ are nonnegative. About properties of positive LSS, see~\cite{FM, LA} and references therein.
To characterize positive uniform LSS we apply Theorem~\ref{th40} on nonnegative c.s.r. semigroups.
\begin{theorem}\label{th90}
 A positively irreducible LSS is uniform with~$\sigma = 0$
 if and only if the family $I + \cA$ has a proper affine invariant subspace~$V \subset \re^d$
 that intersects the positive orthant and does not contain the origin.  In this case,
 every trajectory staring in the set~$\fP = V \cap \re^d_+$ never leaves~$\fP$.
 \end{theorem}
{\tt Proof}. If the LSS is uniform, then the family~$\cM$ is c.s.r. Since this is a positively irreducible
family of nonnegative matrices, Theorem~\ref{th40} implies the existence of an affine
subspace~$V$ such that $0\notin V$, the set $\fP = V \cap \re^d_+$ is nonempty and compact, and
$e^{tA}\fP \subset \fP, \, A \in \cA, \, t \in \re_+$.
The latter means that for every switching law
$A(t)$, if $x(0) \in \fP$, then $x(t) \in \fP$ for all $t > 0$. Hence, $\dot x(t) \in \tilde V$, and therefore
$A(t)x(t) \in \tilde V$ for all~$t > 0$ (let us recall that $\tilde V$ is the linear part of~$V$).
On the other hand, $x(t) \in V$, and so
$(I + A(t)) x(t) \in x(t) + \tilde V = V$. Thus, $V$ is invariant for $I + \cA$, and the trajectory
$x(t), \, t > 0$, lies in~$\fP$. Conversely, if the desired subspace~$V$ exists,
then $A  V \subset \tilde V, \, A \in \cA$, hence,  $\dot x(t) = A(t)x(t) \in \tilde V$, whenever
$x(t) \in V$. Thus, every trajectory~$x(t)$ starting in~$V$ stays in~$V$.
Therefore, $e^{tA}V \subset V$ for every $e^{tA} \in \cM$, and, by Theorem~\ref{th40},
the family~$\cM$ is c.s.r.

 {\hfill $\Box$}
\smallskip

\smallskip
\begin{center}
{\tt 9.3. Regularity of fractal curves}
\end{center}
\smallskip

Let $\{B_0, B_1\}$ be an irreducible pair of affine contractions in~$\re^d$.
The irreducibility means that they do not share a common {\em affine} invariant plane.
Since $B_i$ is a contraction, it has a unique fixed point: $B_iv_i = v_i\, , \, i = 0,1$.
Assume also the {\em Barnsley condition} ({\em cross-condition}): $B_0v_1 = B_1v_0$.
Then the following functional equation:
\begin{equation}\label{fract}
v(t) \ = \
\left\{
\begin{array}{lcl}
B_0 v(2t) &, & t \in \bigl[0, \frac12 \bigr];\\
B_1 v(2t-1) &, & t \in \bigl[\frac12 , 1 \bigr]\,
\end{array}
\right.
\end{equation}
possesses a unique continuous solution $v: [0,1]\to \re^d$ called a {\em fractal curve}~\cite{B, MP}.
We consider the simplest case of equation~(\ref{fract}), with two operators and with double contraction of the argument.
For more general constructions, see~\cite{GS, P2.5, V2}.
The H\"older exponent of the solution is $\alpha_v = -\log_2 \rho(B_0, B_1)$.
The local regularity at each point~$t \in [0,1]$ is measured by the local H\"older exponent:
 $$
 \alpha_v (t) \ = \ \sup \, \Bigl\{\, \alpha > 0 \quad \Bigl| \quad  \bigl|v(t+h) - v(t)\bigr| \, \le \, C \, h^{\, \alpha}\Bigr\}\, .
 $$
By \cite[Theorem 5]{P2.5}, for every point $t \in [0,1]$, we have $\alpha_v (t) \in [\alpha_{\min}, \alpha_{\max}]$ with
$$
\alpha_{\min} \ = \ - \, \log_2 \, \rho \, (B_0, B_1)\ ; \quad  \alpha_{\max} \ = \ -\, \log_2 \, \check \rho \, (B_0, B_1)\, , $$
where, let us remember,
$\rho = \rho_{+\infty}$ is the joint spectral radius and $\rho = \check \rho_{-\infty}$ is the lower spectral radius
(Section~2). Moreover, if $B_0, B_1$ are both nonsingular, then for any $\nu \in [\alpha_{\min}, \alpha_{\max}]$
the set of points $t$ for which $\alpha_v(t) = \nu$ is everywhere dense in the segment~$[0,1]$. This set is of Lebesgue measure
zero for all but one~$\nu$, for which it is of measure~$1$ (this is for $\nu_0 = -\log_2 \rho_0 (B_0, B_1)$,
where $\rho_0$ is the Lyapunov exponent).

Thus, every fractal curve has a varying local regularity. In each iterval $(t_0, t_1) \in [0,1]$,
all values of the local H\"older exponent $\alpha_v(t)$ from $\alpha_{\min}$ to $\alpha_{\max}$  are attained.
The only exception is when $\alpha_{\min} = \alpha_{\max}$, in which case the local regularity is the same at all
points~$t \in [0,1]$. Our goal is to characterize fractal curves with constant local regularity.

The equality~$\alpha_{\min} = \alpha_{\max}$ means that $\check \rho(\tilde B_0, \tilde B_1) = \rho(\tilde B_0, \tilde B_1)$,
which is equivalent, in view of Proposition~\ref{p3}, to the existence of number $r>0$
such that the pair $\{r^{-1}\tilde B_0\, , \, r^{-1}\tilde B_1\}$ is c.s.r.
We use factorization (\ref{blocks}) for matrices of the linear parts $\tilde B_0, \tilde B_1$ of operators~$B_0, B_1$.
Applying Theorem~\ref{th5} and~\ref{th10}, we obtain
\begin{theorem}\label{th120}
A fractal curve generated by nonsingular operators $B_0, B_1$ has a constant local regularity
 if and only if there is a basis in $\re^d$ in which both matrices $\tilde B_0, \tilde B_1$
 have the form~(\ref{blocks}), where in one diagonal block, both matrices $\tilde B_0^{(j)}$ and
 $\tilde B_1^{(j)}$ are orthogonal multiplied by a number~$r\in \bigl[\frac12,1 \bigr)$,
and in the other blocks,  $\rho(B_{0}^{(i)}, B_{1}^{(i)}) \le r, \, i \ne j$.

 In this case $\, \alpha_v(t) = -\log_2 r\, $ at all points $t \in [0,1]$.
\end{theorem}
If the  operators~$B_0, B_1$ have irreducible linear parts, then the situation is much simpler:
\begin{cor}\label{c50}
If $B_0, B_1$ are nonsingular affine operators  whose linear parts do not share a common invariant subspace, then
the fractal curve has a constant local regularity if and only if these operators
  are similarities with the same contraction factor~$r \in \bigl[\frac12,1 \bigr)$.
\end{cor}
\begin{ex}\label{ex40}
{\em The famous Koch curve (``the Koch snowflake'') has constant local regularity $\alpha_v(t) = \log_2 \sqrt{3}$,
because it is generated by two similarities
with $r = 1/\sqrt{3}$. Another famous fractal curve is De Rham curve obtained as a limit of the
cutting angle algorithm from a polygon. In each iteration, all sides of a current polygon are divided into three parts
with the same ratio $\omega : (1-2\omega): \omega$, where $\omega \in (0, 1/2)$. This curve has a constant local regularity
only for $\omega = \frac14$, for other values of~$\omega$ we have $\alpha_{\min} < \alpha_{\max}$ (see~\cite{P2.5}).
Other important examples are refinable functions and wavelets considered in the next subsection.
}
\end{ex}

\smallskip
\begin{center}
{\tt 9.4.  Refinable functions and wavelets}
\end{center}
\smallskip

The following functional equation on a compactly supported scalar function~$\varphi$ plays a crucial role in the construction of compactly
supported wavelets and of subdivision schemes in approximation theory, curve and surface design~(see~\cite{MP, DL, P2.5, V} and references therein):
\begin{equation}\label{eqn.ref}
\varphi(x) \ = \ \sum_{k=0}^{N}\, c_k \, \varphi \, (2t - k)\,  .
\end{equation}
This is called {\em refinemet equation}. Its compactly supported solution~$\varphi$ ({\em refinable function}) is a fixed point
of {\em transition operator} $[Tf](t) = \sum_{k=0}^Nc_k f(2t - k) $. It can be assumed that $\sum_{k=0}^N c_k = 2$, the general case can always be reduced to this one. In this case, equation~(\ref{eqn.ref}) has a unique, up to normalization, solution~$\varphi$ in the space of distributions. We normalize this solution by the condition~$\int_{\re} \varphi(t)d t = 1$. This solution is supported on the segment~$[0, N]$ and for an arbitrary compactly supported distribution~$f$
such that $\int_{\re} f(t) d t = 1$,  one has $T^kf \to \varphi$ as $k \to \infty$ (the convergence is in the sense of distributions). The main problem is when this solution is summable or continuous and what is its regularity and other properties.
A necessary condition for summability of~$\varphi$ are formulated in terms of the generating polynomial
$m(z) = \frac12 \sum_{k=0}^N c_k z^k$.
\smallskip

\noindent \textbf{Theorem A}~\cite{V}. {\em If the solution of refinement equation~~(\ref{eqn.ref}) is summable, then at least one of the following two conditions are fulfilled:} a) {\em $\, m(-1) = 0$}; b) {\em there is $z \in \co$ such that $m(z) = m(-z) = 0$.}
\smallskip

 Continuous refinable function  generate wavelet functions in construction of compactly supported wavelets (for instance, Daubechies wavelets),  the limit function of a subdivision scheme, etc. The specific fractal-like properties of refinable functions
are explained by the fact that equation~(\ref{eqn.ref}) becomes the equation on fractal curve~(\ref{fract})
for the vector-function~$v(t) = v_{\varphi}(t)= \bigl(\varphi(t), \ldots , \varphi(t-N+1) \bigr)^T$ with the operators
$B_0, B_1$ given by $N\times N$ matrices:
\begin{equation}\label{toep}
(B_k)_{ij} \ = \ c_{2i-j+k-1}\ , \qquad i, j = 1, \ldots , N\, ;  \quad k = 0,1
\end{equation}
(if the index $2i-j+k-1$ is negative or exceeds~$N$, then we set $c_{2i-j+k-1} = 0$).
This equation is considered on the affine subspace $V \subset \re^N$,
the smallest by inclusion common invariant affine subspace of the matrices~$B_0, B_1$
passing through the point $v(0)$, which is an eigenvector of $B_0$ corresponding to the
eigenvalue~$1$.

Thus, $v_{\varphi}(t)$ is a fractal curve of operators~$B_0, B_1$ on~$V$. In particular, it has a varying local regularity
   in the sense described in the previous subsection. The same property is inherited by all refinable functions~$\varphi(t)$, wavelet functions,    and limit functions of subdivision schemes. Our aim is to
   characterize refinable functions with the constant local regularity, i.e., the case $\alpha_{\min} = \alpha_{\max}$.
 \begin{conj}\label{conj10}
 If a refinable function has constant local H\"older exponent~$\alpha$, then~${\alpha = 1}$.
 \end{conj}
Thus, either a refinable function $\varphi$ is smooth (for instance, Lipschitz continuous), in which case
$\alpha_{\varphi}(t) = 1$ for all~$t$, or $\varphi$ has varying local regularity and $\alpha_{\min} < \alpha_{\max}$.
We are able to prove Conjecture~\ref{conj10} only under some mild assumption:
\begin{prop}\label{p120}
Conjecture~\ref{conj10} holds true unless all roots of~$m(z)$ are on the unit circle $\{z \in \co \ | \ |z| = 1\}$.
\end{prop}
{\tt Proof. } We use several results from~\cite{P3}. There is a basis in the space~$V$ such that
the matrices $B_i|_{V}$ have block lower triangular form with numbers~$\frac12, \ldots , \frac{1}{2^q}$ and a block
$\hat B_i$ on the diagonal, $i=0,1$. The matrices $\hat B_i$ are both nondegenerate and the pair~$\{\hat B_0, \hat B_1\}$ is irreducible. These matrices have the same form~(\ref{toep}) with a special sequence~$\hat c_1, \ldots , \hat c_{n}$;
the generating polynomial of this sequence $\hat m(z) =\frac12 \sum_{k=0}^n \hat c_k z^k$  have all roots on the
unit circle if and only if the polynomial $m(z)$ does. From the structure of the matrix~$B_i|_{V}$
it follows that if $r = \rho(\hat B_0, \hat B_1) \le \frac12$, then $\check \rho (B_0|_{V}, B_1|_{V}) =
\rho (B_0|_{V}, B_1|_{V}) = \frac12$, hence $\alpha_{\min} = \alpha_{\max} = 1$, which completes the proof.
If $r > \frac12$, then $\alpha_{\min} = \alpha_{\max}$ if and only if the pair~$\{r^{-1}\hat B_0, r^{-1}\hat B_1\}$
is c.s.r. Since this pair is irreducible and the matrices are nonsingular, it follows from Theorem~\ref{th10} that
both these matrices are orthogonal in some basis. Hence, the matrix $\hat B_0^{-1}\hat B_1$ is orthogonal in that basis.
In particular, all its eigenvalues are on the unit circle. On the other hand, the eigenvalues of this matrix
are minus roots of the polynomial~$\hat m(z)$~\cite[Lemmas 1,2]{P3}. Thus, all roots of $\hat m(z)$ and hence all
roots of~$m$ are on the unit circle.

 {\hfill $\Box$}
\smallskip

\smallskip
\begin{center}
{\tt 9.5. The Euler binary partition function}
\end{center}
\smallskip

For arbitrary  $r \in \n \cup \{\infty \}$, the Euler binary partition function
$b(k) = b_r(k)$ is defined on the set of nonnegative integers
$k$ as the total number of different binary expansions $k =
\sum_{j = 0}^{\infty } d_j 2^j$, where the ``digits'' $d_j$ take
values from the set $D = \{ d \in \z \ | \ 0\le  d < r\}$.
For instance, if $r = 2$, then $D = \{0,1\}$, and every number $k$ has a unique
binary expansion, therefore $b_2(k) \equiv 1$. If $r \ge 3$, then a binary expansion may not be unique,
and a natural question arises how does $b_r(k)$ grow with $k$~?
 The asymptotic
behavior of $b_r(k)$ as $k \to \infty$ was studied in various
contexts by L.~Euler, K.~Mahler, N.G.~de Bruijn, D.E.~Knuth,
 and others (see~\cite{R} for  bibliography and historical comments).
 Leonard Euler first considered the function $b_{\infty}( k)$ in connection with
 the power series $\sum_{k=0}^{\infty} b_{\infty}(k) z^k \, = \, \prod_{j=0}^{\infty}\frac{1}{1 - z^{2^j}}$.
 In 1940 K.~Mahler proved that $b_{\infty} (k) \, \asymp \, k^{\, \frac12 \log_2 k}$, this asymptotic
 relation was sharpened by N.G.~de Bruijn (1948) and by  D.E.~Knuth (1966). The case of finite~$r$
 was first systematically analyzed by~B.Reznick in 1990~\cite{R}. He showed that for every even~$r$,
 we have $b_r(k)\, \asymp \, k^{\, \log_2 (r/2)}$, i.e., the partition function grows polynomially in~$k$.
  For odd~$r$, however, the asymptotics of $b_r(k)$ as $k \to \infty$ is not that regular. Denote
\begin{equation}\label{p1p2}
p_1(r) \quad =\quad \liminf_{k \to \infty }\ \log b_r(k)/ \log k;
\qquad p_2(r) \quad = \quad \limsup_{k \to \infty }\ \log b_r(k)/ \log
k\, .
\end{equation}
Thus, for even $r$ we have $p_1 =p_2 = \log_2 (r/2)$.
However, already for $r=3$, the lower and upper exponents of growth do not coincide $p_1(3)=0$ and $p_2(3) = \log_2 \frac{\sqrt{5}+1}{2}$.
B.Reznick in~\cite{R} formulated the following problem: is it true that $p_1 < p_2$ for all even $r$~?
Applying Theorem~\ref{th40} we obtain  an affirmative answer:
\begin{theorem}\label{th100}
For every odd $r$ we have $p_1(r) < p_2(r)$.
\end{theorem}
Thus, for every odd $r$ the lower and upper exponents of growth of the function $b_r(k)$
do not coincide. To prove the theorem we need several auxiliary results and observations.
In~\cite{P2} it was shown that $p_{1} \, = \, \log_2
\check \rho(D_0, D_1)\, $ and $\, p_{2} \, = \, \log_2  \rho(D_0,
D_1)$, where $D_0, D_1$ are $(r-1)\times (r-1)$-matrices defined
as follows:
\begin{equation}\label{2.T}
(D_s)_{ij} \quad  = \quad
\left\{
\begin{array}{ccl}
1 \, , \, & \, \mbox{if} \, & \, 1-s \, \le \, 2i-j \, \le \, r-s\, , \\
0\, , \, & \, \mbox{otherwise}
\end{array}
\right.
\end{equation}
(here $s = 0,1$ and $i, j \in \{1, \ldots , r-1\}$).
For example, for $r=5$ we have the following $4\times 4$-matrices:
$$
D_0\quad = \quad
\left(
\begin{array}{rrrr}
1 & 0& 0& 0\\
1 & 1& 1& 0\\
1 & 1& 1& 1\\
0 & 0& 1& 1
\end{array}
\right)\ ; \qquad
D_1\quad = \quad
\left(
\begin{array}{rrrr}
1 & 1& 0& 0\\
1 & 1& 1& 1\\
0 & 1& 1& 1\\
0 & 0& 0& 1
\end{array}
\right)
$$
Consider the matrices $B_s \, = \, \frac2r \, D_s, \, s = 0,1$. Note that the matrix $\bar B = \frac12 \bigl( B_0+B_1\bigr)$
is column-stochastic, and hence $\rho_1 (B_0,B_1) \, = \, \rho  \bigl( \bar B \, \bigr) \, = \, 1$.
Therefore, $\, \check \rho (B_0,B_1)\, \le \, 1 \, \le \, \rho (B_0,B_1)$. We see that
 $p_1 = p_2$ if and only if $\check \rho \, = \,  \rho \, = \, 1$.
 Thus, we have proved
\begin{prop}\label{p100}
We have $p_1(r) = p_2(r)$ if and only if the family $\{B_0, B_1\}$ is c.s.r.
\end{prop}
From this fact it is already clear why $p_1(r) = p_2(r)$ for even~$r$. Indeed,
each column of the matrices~$D_0, D_1$,   has exactly $\frac{r}{2}$ ones, hence
$B_0, B_1$ are both column-stochastic, and the family $\{B_0, B_1\}$ is c.s.r. Theorem~\ref{th100}
states that that this family is not c.s.r. for odd~$r$.
\smallskip

{\tt Proof of Theorem~\ref{th100}}. Comparing formulas~(\ref{2.T}) and~(\ref{toep}) we see that our matrices $B_0, B_1$ correspond to the refinement equation~$T\varphi = \varphi$ with the transition operator $[T\varphi ](t) \, = \, \frac2r \, \sum_{k=0}^{r-1}\, \varphi (2t-k)$. The associated generating polynomial is~$\, m(z) \, = \, \frac{1}{r}\, \sum_{k=0}^{r-1}z^k \, = \, \frac{z^r - 1}{r(z-1)}$.
If $r$ is odd, then $m(-1) = \frac1r \ne 0$. Moreover, $\bigl(z^r - 1\bigr) \, + \, \bigl(\, (-z)^r -1 \bigr)\, \equiv \, -2$, hence the values $m(z)$ and $m(-z)$ cannot vanish simultaneously. Thus, by Theorem~A, the refinement equation does not have a nontrivial summable solution.

Assume now that $\{B_0, B_1\}$ is c.s.r. The matrix~$\bar B$ has strictly positive elements on the main diagonal and on both the neighboring diagonals, hence $\bar B$ is  positively irreducible, and so
is the pair~$\{B_0, B_1\}$. Theorem~\ref{th40} implies that $B_0$ and $B_1$ possess a
a common invariant  affine subspace $V \subset \re^{r-1}$ such that the set  $\, \fP \, = \, V \cap \re^{r-1}$ is bounded.
Since $B_i\fP \subset \fP$, it follows that the transition operator respects the set of functions
$f \in L_1(\re)$ such that $v(t) = \bigl(\varphi(t), \ldots , \varphi(t-r+2) \bigr)^T \in \fP$
for almost all~$t \in [0,1]$.
Taking an arbitrary initial function~$f_0$ from this set,
we see that $T^nf_0$ belongs to this set for each $n$, and hence, so does the limit $\varphi = \lim_{n\to \infty}\, T^nf_0\, $,
which is, by Theorem~A, a solution of the refinement equation. For this solution $v(t) \in \fP$ a.e., therefore $\varphi \in L_{\infty}[0, 1]$. This contradiction completes the proof.

{\hfill $\Box$}

\smallskip

\bigskip

\begin{center}
\textbf{10. Appendix}
\end{center}
\medskip

{\tt Proof of Proposition~\ref{p5}}. The necessity is obvious, let us prove sufficiency. If the spectral radii
of matrices from~$\cS$ are bounded, then they are bounded by one. Assume to the
contrary that there are $A_i \in \cS$ such that $\|A_i\| \to \infty$ as $i \to \infty$.
Since~$\cS$ has a finite irreducible subset, we assume that the sequence
$\{A_i\}_{i=1}^{\infty}$ starts with this subset, i.e., the set $\{A_1, \ldots , A_n\}$
is irreducible.  For each $k \ge \n$ we denote by $\cS_k$ the semigroup generated by the matrices $\{A_1, \ldots , A_k\}$.
This set is irreducible, and since the spectral radii of all products of these matrices do not exceed one, their
joint spectral radius is at most one. Hence, by Barabanov's theorem~\cite{B1},
 there is a norm $\|\cdot\|_k$ in $\re^d$ such that
$\|A_i\|_k \le 1,  \, i = 1, \ldots , k$. We normalize this norm by the condition that its maximum on the
unit Euclidean ball is one. By the compactness, there is a limit point $f(\cdot )$ of this
sequence of norms. Passing to a subsequence, it may be assumed that $\|x\|_k \to f(x)\, , \ x \in \re^d$,
 as $k \to \infty$. Obviously, $f$ is a nonnegative  convex symmetric homogeneous functional.
 The set of points $x \in \re^d$ such that $f(x) = 0$ is a proper
 linear subspace of $\re^d$ invariant with respect to all $A_i$, which contradicts the irreducibility.
 Hence, $\|x\| > 0$ for $x\ne 0$, and therefore $f$ is a norm. The induced operator norm
 of all $A_i$ is at most one. Hence, by the equivalence of norms, $\|A_i\|$ cannot tend to infinity as $i \to \infty$.

{\hfill $\Box$}
\medskip

\medskip

{\tt Proof of Proposition~\ref{p40}.} If $\cS$ is positively irreducible, then it contains a finite positively irreducible subset.
Therefore, without loss of generality it can be assumed that $\cS$ is generated by a finite matrix family $\cA = \{A_1, \ldots , A_m\}$.
The matrix $\bar A = \frac1m \sum_{i=1}^m A_m$.
 is positively irreducible, and $\bar A V \subset V$, hence $V$ intersects the interior of~$\re^d_+$.
 Otherwise, the set~$\fP$ is contained in some coordinate subspace
 which is invariant for~$\bar A$. If~$\fP$ is unbounded, then it contains a ray. A parallel ray starting at the origin
  is contained in the subspace~$\tilde V$ (the linear part of~$V$). Consequently, the cone~$K = \re^d_+\cap \tilde V $
   is nontrivial. Since~$\bar A K \subset K$, it follows from the Krein-Rutman theorem~\cite{KR}
   that~$\bar A$ possesses an eigenvector~$v \in K$.
 Hence, $v$
is a Perron-Frobenius eigenvector of the matrix~$\bar A$, because a positively irreducible matrix has a unique, up to normalization, nonnegative eigenvector. Thus, $V$ is parallel to the vector~$v$. Let $\bar A^T$ be the transpose matrix and $v^*$ be its Perron-Frobenius eigenvector. For sufficiently large $s > 0$, the set $\, V_s\, = \, \{x \in \re^d_+ \ | \
 (v^*, x) = s\}$ intersects $V$ by some nonempty set $\fG$. Clearly, $\bar A \, \fG \, \subset \, \fG$. Since $\bar A^T$ is positively irreducible, it follows that~$v^* > 0$,
  and hence the set $\fG$ is bounded. The operator $\bar A$ respects the cone~$K_s = \{tx \ | \ t \ge 0\, , \ x \in \fG\}$, and consequently, $v \in K_s$. Hence, $\tau v \in V$ for some $\tau > 0$. However, $\tau \, v \in \tilde V$, and so
  $0 = (\tau v \, - \, \tau v) \, \in \, V$, which contradicts the assumption.

{\hfill $\Box$}

\smallskip

 {\tt Proof of Theorem~\ref{th40}}. Sufficiency. Suppose~$\cS$ has a common invariant affine subspace~$V$
satisfying all the assumptions. By Proposition~\ref{p40}, the set~$\fP \, = \, V \cap \re^d_+$ is bounded
and contains some interior point  $z \in {\rm int}\, (\re^d_+)$. There are positive constants $C_1, C_2$
such that $\, C_1 \|Bz\| \, \le \, \|B\|\, \le \, C_2 \|Bz\|$ for each nonnegative matrix~$B$.
Since the set $\{Az\ | \  A \in \cS\}$ is contained in~$\fP$,
we see that $\, C_1 \inf\limits_{x \in \fP}\|x\| \, \le \, \|A\|\, \le \, C_2 \sup\limits_{x \in \fP}\|x\|$
for all $A \in \cS$. Thus,~$\cS$ is bounded and separated from zero, hence it is c.r.s.

Necessity. Since~$\cS$ is irreducible, it contains a finite positively irreducible subset~$\cB$. Every finite subset~$\cA \subset \cS$ that contains~$\cB$ is also irreducible. If every such a subset~$\cA$ possesses a common invariant affine subspace satisfying all the assumptions, then so does the set~$\cS$. Thus, it suffices to prove the necessity for finitely generated c.s.r. semigroups. Assume $\cS = \cS_{\cA}$ for some~$\cA = \{A_1, \ldots , A_m\}$.
For a nonnegative matrix family~$\cA$ we have $\rho_1(\cA) \, = \, \rho(\bar A)$, where $\bar A = \frac1m \sum_{i=1}^m A_m$ (see~\cite{P4}). Hence, in our case,  $\rho(\bar A) = 1$. Furthermore, since $\cA$ is positively irreducible,
it follows that there is a norm~$\|\cdot \|$ such that $\max\limits_{i=1, \ldots , m}\|A_ix\| \, = \, \rho(\cA)\, \|x\|$,
for all $x  \in \re^d_+$ (see~\cite[Theorem 3]{GP}).  In our case~$\rho(\cA) = 1$, hence
\begin{equation}\label{sum0}
\max_{d_1, \ldots , d_k}\|A_{d_k}\cdots A_{d_1}x\| \ = \ \|x\|\ , \qquad x \in \re^d_+\, , \, k \in \n\, .
\end{equation}
Let~$v \in {\rm int}(\re^d_+)\, , \, \|v\| = 1$,
be a Perron-Frobenius eigenvector of~$\bar A$. We have $v = \bar A v$, therefore
\begin{equation}\label{sum1}
v \ = \ \bar A^kv \ = \ m^{-k}\sum_{d_1, \ldots , d_k}A_{d_k}\cdots A_{d_1}v\, , \qquad k \in \n\, .
\end{equation}
Let $L_v \, = \,
\{x \in \re^d \ | \ (p, x) = 1\}$ be  the hyperplane of support for the unit ball of the norm~$\|\cdot\|$
passing through the point~$v$. Thus, $(p, v) = 1$ and $(p, x) \le 1$ for every $x \in \re^d, \, \|x\| \le 1$.
 Substituting in~(\ref{sum1}), we get
\begin{equation}\label{sum2}
 1 \ = \ (p, v)\ = \ m^{-k}\sum_{d_1, \ldots , d_k}\, \bigl(\, p \, , \, A_{d_k}\cdots A_{d_1}v\, \bigr) \, .
\end{equation}
By~(\ref{sum0}) we have $\|A_{d_k}\cdots A_{d_1}v\|\, \le \, \|v\| \, = \, 1$, therefore,
$\bigl(\, p \, , \, A_{d_k}\cdots A_{d_1}v\, \bigr) \le 1$, for every $d_k, \ldots , d_1$.
Combining with~(\ref{sum2}) we obtain $\, ( p \, , \, A_{d_k}\cdots A_{d_1}v)\, = \, 1$, for every $d_k, \ldots , d_1$,
which means that all the points $A_{d_k}\cdots A_{d_1}v$ belong to the affine hyperplane~$L_v$.
Let $V$ be the affine hull of the points $v$ and $A_{d_k}\cdots A_{d_1}v$ over all possible
 products of lengths $k \in \n$. Clearly, $A_iV \subset V\, , \ i = 1,\ldots , m$
and $v \in V \subset L_v$. Thus, $V$ is a required proper affine subspace.

{\hfill $\Box$}

\smallskip

{\tt Proof of Lemma~\ref{l10}.} Let $C \in \cS$ be an operator of the minimal rank~$r$ and let $L$ be the range of~$C$.
Clearly, $C$ is nonsingular on~$L$, otherwise the rank of~$C^2$ is smaller than~$r$. Hence ${\rm dim }\, L \, = \,
r$, the kernel $K$ of $C$ is of dimension~$d-r$ and does not intersect~$L$.
 Since $C$ is a direct sum of a nilpotent operator and of an operator similar to orthogonal~\cite[Theorem 2.5]{OR}, the operator
 $C|_{L}$ is similar to orthogonal. Hence, powers of this operator come arbitrarily close to the identity
 operator. Therefore, some powers $C^n$ are arbitrarily close to the projector~$P$ onto $L$ parallel to
 the subspace~$K$. Since $\cS$ is closed, it contains~$P$. Consider the semigroup $\cS_{P} = \bigl\{(PA)|_{L}\ \bigl| \
 A \in \cS\bigr\}$.  The range of the operator $(PA)|_{L}$ coincides with the range of $PAP$, whose dimension is~$r$,
 otherwise ${\rm rank}(PAP) < r$. Whence, $(PA)|_{L}$ is nonsingular. Moreover, $\rho\bigl((PA)|_{L}\bigr) \, = \,
 \rho\bigl(PAP\bigr)\, = \, 1$.
 Thus, the semigroup~$\cS_P$ consists of nonsingular matrices and has constant spectral radius. Let us show that $\cS_P$ is irreducible.
 If, to the contrary,  there is a nontrivial subspace $L' \subset L$ such that $PA (L') \subset L'$
 for all $A \in \cS$, then $AL' \subset (L' + K)$. Hence, the linear span of the set
 $\bigl\{ AL' \ \bigl| \ A \in \cS\bigr\}$, which is a common invariant subspace for all elements of~$\cS$,
 is contained in the subspace~$L'+ K$. This contradicts the irreducibility of~$\cS$. Thus,
 $\cS_P$ is irreducible, is of constant spectral radius, and consists of nonsingular operators.
 By Theorem~\ref{th10}, one can choose a basis in~$L$ so that all operators of~$\cS_P$ become orthogonal.
      Now we complement this basis to a basis of~$\re^d$ so that $K$ becomes orthogonal to~$L$.
 In the new basis,~$P$ is an orthogonal projector, which concludes the proof.

{\hfill $\Box$}

\end{document}